\documentclass[10pt]{article}

\usepackage{epsfig}
\usepackage{amsmath}
\usepackage{amsfonts}
\usepackage{amsthm}
\usepackage{amsbsy}
\usepackage{appendix}
\usepackage{epstopdf}
\usepackage[table]{xcolor}
\usepackage{subfloat}
\usepackage{float}
\usepackage[thinlines]{easytable}
\usepackage{array}
\usepackage{color}

 \allowdisplaybreaks

        \theoremstyle{plain}

        \theoremstyle{remark}
        \newtheorem{remark}{\bf Remark}[section]
        \theoremstyle{remark}
        
        \theoremstyle{remark}

	\newcommand{\Rho}{\mathrm{P}}

\title{Efficient Stochastic Asymptotic-Preserving IMEX Methods for Transport Equations with Diffusive Scalings and Random Inputs\footnote{This research was supported by NSF grants DMS-1522184 and
DMS-1107291: RNMS KI-Net, by NSFC grant No. 91330203, and by the Office of the Vice
Chancellor for Research and Graduate Education at the University of Wisconsin, Madison,
with funding from the Wisconsin Alumni Research Foundation. The research was partially supported by the research grant GNCS-INDAM 2016 {\it Numerical methods for uncertainty quantification in hyperbolic and kinetic equations}.}}
\author{Shi Jin\footnote{Department of Mathematics, University of Wisconsin, Madison, WI 53706, USA (sjin{@}wisc.edu) and Institute of Natural Sciences, School of Mathematical Science, MOE-LSEC and SHL-MAC, Shanghai Jiao Tong
University, Shanghai 200240, China.} , Hanqing Lu\footnote{Department of Mathematics, University of Wisconsin, Madison, WI 53706, USA (hanqing{@}math.wisc.edu).} and Lorenzo Pareschi\footnote{Department of Mathematics \& Computer Science, University of Ferrara, Ferrara 44121, Italy (lorenzo.pareschi{@}unife.it).}} 

\begin{document}
\maketitle
\bigskip
\begin{center}
{\bf Abstract}
\end{center}

For linear transport and radiative heat transfer equations with random inputs, we develop new generalized polynomial chaos based Asymptotic-Preserving stochastic Galerkin schemes that allow efficient computation for the problems that contain both uncertainties and multiple scales. Compared with previous methods for these problems, our new method use the implicit-explicit (IMEX) time discretization to gain higher order accuracy, and by using a modified diffusion operator based penalty method, a more relaxed stability condition--a hyperbolic, rather than parabolic, CFL stability condition, is achieved in the case of small mean free path in the diffusive regime. The stochastic Asymptotic-Preserving property of these methods will be shown asymptotically, and demonstrated numerically, along with computational  cost comparison with previous methods.
\\



\tableofcontents

\section{Introduction}

Transport equations, for example, linear transport equations and radiative transfer equations, arise in many important physical applications, from the microscopic neutron transport to large scale astrophysical problems \cite{Cha}. As in general kinetic theory, one of the main challenges is its high dimensionality, since these equations describe particle density distributions  defined in the phase space, with independent variables time, space and particle velocity. Another challenge is the multiscale nature, since the Knudsen number, the dimensionless mean free path, could vary in different orders of magnitude in the computational problems. When the Knudsen number is $O(1)$, we are in the kinetic regime.  When it is small, usually one can approximate the kinetic equations by hydrodynamic or diffusion type equations which are macroscopic
equations defined in physical space thus are much more efficient for numerical computations. Asymptotic-preserving (AP) schemes, which mimic the asymptotic transition from the kinetic to the macroscopic equation at the discrete level, have been proven to be an effective computational paradigm for such multiscale kinetic or transport equations, and have found many applications in kinetic theory \cite{Deg1, Deg2, Acta,Jin_rev}. For those relevant to the equations to be studied in this article, see for example in neutron transport equations \cite{boscarino2013implicit,CGLV,GT,jin2000uniformly, Klar, LS, LMM1, LMM2, LM, NP}, and  nonlinear  radiative transfer \cite{BC, klar2001numerical, JiangXu}.

In this article, we are interested in addressing the third challenge in kinetic modeling: the {\it uncertainty}, for transport and radiative transfer equations.  Since kinetic equations often arise from microscopic equations, for example Newton's second law or particle dynamics, via taking  mean-field limits \cite{BGP, Spohn},  one often encounters difficulties in determining precisely or accurately terms such as collision kernels or scattering coefficients.  In addition, initial, boundary data or forcing or source terms could
also be measured inaccurately. Thus uncertainty is {\it intrinsic} in kinetic models that must be studied in order to assess the accuracy of kinetic models and validate and enhance the computational reliability of kinetic modelings.

While uncertainty quantification (UQ) has been a popular field of scientific computing in many areas of scientific and engineering research, its study for kinetic equations has been scarce. Only recently one sees some  development of efficient numerical methods for kinetic equations with uncertainties \cite{ChenLiuMu, HJ-Boltzmann, JinLiu, JinLu,  jin2015asymptotic, ZhuJin}.  In particular, in \cite{jin2015asymptotic},
the notion of {\it stochastic Asymptotic-Preserving} (sAP) was introduced. For  generalized polynomial chaos (gPC) based stochastic Galerkin (SG) methods, which are the methods to be used in this article to study the uncertainties, it requires an SG method for kinetic equations with uncertainty, as Knudsen number goes to zero, becomes a  SG method for the macroscopic hydrodynamic or diffusion equations. This is the stochastic extension of the deterministic AP schemes. Moreover, as realized in \cite{jin2015asymptotic},  and subsequently in other kinetic equations \cite{JinLiu, JinLu, ZhuJin}, under the SG framework, one obtains a vectorized system of deterministic kinetic equations based on which the deterministic AP machinery can be easily utilized to develop sAP schemes for random kinetic equations. The gPC approach for uncertain kinetic
equations is natural, and efficient, since kinetic equations typically preserve the regularity of the initial data in the random space \cite{HJ-Boltzmann, JLM, JinLiu, LiWang}, thus one indeed obtains {\it spectral} accuracy in the random space, and then the  sAP property allows one to use {\it Knudsen number independent} time
step, mesh sizes and the order of polynomial degrees in the gPC method.

While sAP schemes have been developed for linear transport equation \cite{JLM, jin2015asymptotic} and radiative heat transfer equations \cite{JinLu} with uncertainties, in this article we aim to {\it improve the accuracy and efficiency} of these methods. Since the gPC-SG  based numerical methods lead to much larger system of equations (compared to the original kinetic equations which are scalar equations), and the random variables that describe the uncertainties often live in high dimension,  more efficient numerical methods can help to significantly reduce the computational cost. To this aim, we will make two  improvements. First, we utilize the implicit-explicit (IMEX) Ruge-Kutta (RK) time marching methods (see \cite{boscarino2013implicit, boscarino2016, A} and the references therein) in order to gain high order of accuracy in time discretization. Second, since these equations, due to high scattering rates, encounter diffusive regimes, a typical AP scheme based on explicit discretization of the convection term needs a parabolic type CFL stability condition $\Delta t=O(\Delta x^2)$, where $\Delta t$ is the time step and $\Delta x$ is the spatial mesh size. A fully implicit scheme will overcome this constraint, but will introduce algebraic difficulties when inverting the entire transport operator numerically. We use a diffusion operator penalized method introduced in \cite{boscarino2013implicit}, which allows to improve the CFL condition to a hyperbolic one $\Delta t=O(\Delta x)$. We modify the approach of \cite{boscarino2013implicit}
by replacing the microscopic probability density distribution in the penalty term by the macroscopic density, which further simplifies the inversion of the implicit terms.

The rest of the paper is organized as follows.
We develop the gPC-SG based IMEX-RK sAP schemes for the random linear transport equation in Section 2 and the random radiative heat transfer equations in Section 3. Numerical tests can be found in Section 4 that verify the efficiency and sAP property of these methods. Some final considerations are reported in Section 5.

\section{The Linear Neutron Transport Equation with Random Inputs}
In this section we consider a one-dimensional linear neutron transport equation with random inputs. The randomness may come from the cross section, initial data and boundary data. Let $f(t,x,v,z)$ be the probability density distribution of particles at position $x\in D\subset \mathbb R$, time $t\in \mathbb R^+$ and depending on $v\in(-1,1)$, which is the cosine of the angle between the particle velocity and position. $z\in I_z\subset \mathbb{R}^d (d\geq 1)$ is the random variable with compact support $I_z$. Then the one-dimensional transport equation is 
\begin{equation}\label{1}
\varepsilon\partial_t f+v\partial_x f=\frac{1}{\varepsilon}\left(\frac{\sigma_s}{2}\int_{-1}^1f(v')dv'-\sigma f\right)+\varepsilon Q,
\end{equation}
where $\sigma=\sigma(x,z)$ is the total cross section and $\sigma_s=\sigma_s(x,z)$ is the scattering cross-section. Then the absorption coefficient $\sigma_A=\sigma_A(x,z)$ is defined by $\sigma_A=\frac{\sigma-\sigma_s}{\varepsilon^2}$. $Q=Q(x,z)$ is the source term and the small parameter $\varepsilon$ is the mean free path. Without loss of generality, $\sigma_A$ and $Q$ are neglected in this study.

The well-known diffusion limit $\varepsilon\to 0$ of the one-dimensional deterministic transport equation is of the following form as \cite{LK, BSS}:
\begin{equation}\label{2}
\partial_t \rho=\partial_x\left(\frac{1}{3\sigma_s}\partial_x\rho\right),
\end{equation}
where 
$$\rho=\frac{1}{2}\int_{-1}^1fdv.$$

One can understand the diffusion limit of the transport equation with random inputs through the following even-odd decomposition \cite{jin2000uniformly}. First write equation (\ref{1}) for $v>0$:
\begin{subequations}\label{3}
\begin{align}
&\varepsilon \partial_t f(v)+v\partial_x f(v)=\frac{1}{\varepsilon}\left(\frac{\sigma_s}{2}\int_{-1}^1f(v')dv'-\sigma f(v)\right),\\
&\varepsilon \partial_t f(-v)-v\partial_x f(-v)=\frac{1}{\varepsilon}\left(\frac{\sigma_s}{2}\int_{-1}^1f(v')dv'-\sigma f(-v)\right).
\end{align}
\end{subequations}

Now denote the even and odd parities
\begin{subequations}\label{4}
\begin{align}
&r(t,x,v,z)=\frac{1}{2}(f(t,x,v,z)+f(t,x,-v,z)),\\
&j(t,x,v,z)=\frac{1}{2\varepsilon}(f(t,x,v,z)-f(t,x,-v,z)).
\end{align}
\end{subequations}

Then (\ref{3}) becomes
\begin{subequations}\label{5}
\begin{align}
&\partial_t r+v\partial_x j=-\frac{\sigma_s}{\varepsilon^2}(r-\rho)\\
&\partial_tj+\frac{1}{\varepsilon^2}v\partial_x r=-\frac{1}{\varepsilon^2}\sigma_sj.
\end{align}
\end{subequations}

As $\varepsilon\to 0^+$, (\ref{5}) yields
\begin{equation}\label{6}
r=\rho, \ \ \ j=-\frac{v}{\sigma_s(x,z)}\partial_x r.
\end{equation}

Substituting (\ref{6}) into the first equation of (\ref{5}) and integrating over $v\in [0,1]$, one gets the limiting diffusion equation with random inputs (\ref{2}).

Based on this even-odd decomposition, Jin, Xiu and Zhu proposed a gPC-SG formulation for the linear transport equation (\ref{1}) with random cross-section $\sigma_s(x,z)$. The  Stochastic Asymptotic-Preserving (sAP) property is also 
introduced for this problem and shown both asymptotically and numerically 
for this framework \cite{jin2015asymptotic}. See also \cite{JLM} for a sharp
regularity and sAP proof. However, a parabolic CFL condition $\Delta t=O(\Delta x^2)$ has to be satisfied for the fully discretized scheme in (\ref{5})
since the scheme  will become an {\it explicit} scheme for the limiting diffusion equation (\ref{2}) as $\varepsilon$ goes to $0$.

Now we propose an efficient sAP method for system (\ref{5}) by applying the IMEX Runge-Kutta scheme \cite{boscarino2013implicit} to the gPC Galerkin formulation to get rid of the parabolic CFL condition. Our scheme will fasten the algorithm significantly especially for the vectorized system after the Galerkin projection and can be shown to be sAP.

\subsection{The gPC-SG Formulation}
To approximate the solution, we use the gPC expansion via an orthogonal polynomials series. That is, for random variable $z\in \mathbb R^d$, one seeks
\begin{subequations}\label{7}
\begin{align}
&r(t,x,v,z)\approx r_N(t,x,v,z)=\sum_{k=1}^K\hat{r}_k(t,x,v)\Phi_k(z),\\
&j(t,x,v,z)\approx j_N(t,x,v,z)=\sum_{k=1}^K\hat{j}_k(t,x,v)\Phi_k(z),
\end{align}
\end{subequations}
where $\left\{\Phi_k(z),1\leq k\leq K, K=\begin{pmatrix}
d+N\\
d
\end{pmatrix}\right\}$ are from $\mathbb P_N^d$, the $d$-variate orthogonal polynomials of degree up to $N\geq 1$, and orthonormal
\begin{equation}\label{8}
\int_{I_z}\Phi_i(z)\Phi_j(z)\pi(z)dz=\delta_{ij}, \ 1\leq i,j\leq K=dim(\mathbb P_N^d).
\end{equation}
Here $\pi(z)$ is the probability density function of $z$ and $\delta_{i,j}$ the Kronecker delta function.

Now one inserts the approximation (\ref{7}) into the governing equation (\ref{5}) and enforces the residue to be orthogonal to the polynomial space spanned by$\{\Phi_1,\cdots, \Phi_K\}$. Then, one obtains a set of vectorized deterministic equations for $\hat{\bold r}=(\hat r_1,\cdots, \hat r_K)^T$ and $\hat{\bold j}=(\hat j_1,\cdots, \hat j_K)^T$:
\begin{subequations}\label{9}
\begin{align}
&\partial_t \hat{\bold r}+v\partial_x \hat{\bold j}=-\frac{1}{\varepsilon^2}\bold S(x)(\hat{\bold r}-\hat{\boldsymbol \rho}),\\
&\partial_t\hat{\bold j}+\frac{v}{\varepsilon^2}\partial_x \hat{\bold r}=-\frac{1}{\varepsilon^2}\bold S(x)\hat{\bold j},
\end{align}
\end{subequations}
where 
$$\hat{\boldsymbol \rho}(x,t)=\langle \hat{\bold r}\rangle=\int_0^1\hat{\bold r}dv,$$
and $\bold S(x)=(s_{ij}(x))_{1\leq i,j\leq K}$ is $K\times K$ symmetric and positive definite matrix with entries
\begin{equation}\label{10}
s_{ij}(x)=\int_{I_z}\sigma_s(x,z)\Phi_i(z)\Phi_j(z)\pi(z)dz.
\end{equation}

As $\varepsilon\to 0^+$,
\begin{equation}\label{11}
\hat{\bold r}=\hat{\boldsymbol \rho}, \ \ \hat{\bold j}=-v\bold S^{-1}\partial_x\hat{\bold r}.
\end{equation}

Inserting above into (\ref{9}a) and integrating over $v$, one can obtain:
\begin{equation}\label{12}
\partial_t \hat{\boldsymbol \rho}=\frac{1}{3}\partial_x[\bold S^{-1}(x)\partial_x\hat{\boldsymbol \rho}].
\end{equation}

\begin{remark}
If one applies the gPC-SG formulation for the limiting diffusion equation (2) directly, one gets
\begin{equation}\label{13}
\partial_t\hat{\boldsymbol \rho}=\frac{1}{3}\partial_x[\tilde{\bold S}(x)\partial_x\hat{\boldsymbol \rho}],
\end{equation}
where $\tilde{\bold S}(x)=(\tilde s_{ij}(x))_{1\leq i,j\leq K}$ is a $K\times K$ symmetric and positive definite matrix with entries
\begin{equation}\label{14}
\tilde s_{ij}(x)=\int_{I_z}\frac{1}{\sigma_s(x,z)}\Phi_i(z)\Phi_j(z)\pi(z)dz.
\end{equation}
Though $\tilde{\bold S}$ and $\bold S^{-1}$ are different matrix, proof in \cite{ShuHuJin} shows that $\tilde{\bold S}\partial_x\hat{\boldsymbol \rho}$ and $\bold S^{-1}\partial_x\hat{\boldsymbol \rho}$ are of spectral error.
\end{remark}

\subsection{An Efficient IMEX Runge-Kutta Scheme}
Following \cite{boscarino2013implicit}, we consider a penalization approach to avoid the parabolic CFL condition of a standard AP scheme in the diffusion limit. By adding and subtracting the term $\frac{\mu}{3}\partial_x(\bold S^{-1}\partial_x \hat{\boldsymbol \rho})$ in (\ref{9}a) and the term $\phi v\partial_x \hat{\bold r}$ in (\ref{9}b), we reformulate the problem into an equivalent form:
\begin{subequations}\label{15}
\begin{align}
\partial_t \hat{\bold r}&=\left[-v\partial_x \hat{\bold j}-\frac{\mu}{3}\partial_x(\bold S^{-1}\partial_x\hat{\boldsymbol \rho})\right]+\left[-\frac{1}{\varepsilon^2}\bold S(\hat{\bold r}-\hat{\boldsymbol \rho})+\frac{\mu}{3}\partial_x(\bold S^{-1}\partial_x \hat{\boldsymbol \rho})\right]\nonumber\\
&=f_1(\hat{\bold r},\hat{\bold j})+f_2(\hat{\bold r}),\\
\partial_t\hat{\bold j}&=-\phi v\partial_x\hat{\bold r}-\frac{1}{\varepsilon^2}\left(\bold S\hat{\bold j}+(1-\varepsilon^2 \phi)v\partial_x\hat{\bold r}\right)=g_1(\hat{\bold r})+g_2(\hat{\bold r},\hat{\bold j}),
\end{align}
\end{subequations}
where
\begin{subequations}\label{16}
\begin{align}
&f_1(\hat{\bold r},\hat{\bold j})=-v\partial_x \hat{\bold j}-\frac{\mu}{3}\partial_x(\bold S^{-1}\partial_x\hat{\boldsymbol \rho}),\\
&f_2(\hat{\bold r})=-\frac{1}{\varepsilon^2}\bold S(x)(\hat{\bold r}-\hat{\boldsymbol \rho})+\frac{\mu}{3}\partial_x(\bold S^{-1}\partial_x\hat{\boldsymbol \rho}),\\
&g_1(\hat{\bold r})=-\phi v\partial_x\hat{\bold r},\\
&g_2(\hat{\bold r},\hat{\bold j})=-\frac{1}{\varepsilon^2}\left(\bold S\hat{\bold j}+(1-\varepsilon^2 \phi)v\partial_x\hat{\bold r}\right).
\end{align}
\end{subequations}
Here we choose $\mu=\mu(\varepsilon)$ such that
\begin{equation}\label{17}
\begin{aligned}
&\lim_{\varepsilon\to0}\mu=1,\\
&\mu= 0 \ \ \ \text{if } \ \varepsilon=O(1);
\end{aligned}
\end{equation}
and $\phi=\phi(\varepsilon)$ such that 
\begin{equation}\label{18}
0\leq\phi\leq\frac{1}{\varepsilon^2}.
\end{equation}

The restriction on $\phi$ guarantees the positivity of $\phi(\varepsilon)$ and $(1-\varepsilon^2\phi(\varepsilon))$ so the problem remains well-posed uniformly in $\varepsilon$. We make the same simple choice of $\phi$ as in \cite{jin2000uniformly}:
\begin{equation}\label{19}
\phi(\varepsilon)=\min\left\{1,\frac{1}{\varepsilon^2}\right\}.
\end{equation}

\begin{remark}
Compared to the penalty term $\mu v^2\partial_x(\bold S^{-1}\partial_x\hat{\bold r})$ used in \cite{boscarino2013implicit}, our choice of the penalty term in (\ref{9}a) is {\it independent} of $v$. So only one large tridiagonal matrix needs to be inverted later at each stage of the IMEX-RK scheme for the fully discretization instead of a series of large tridiagonal matrix (the number depends on the number of $v$ used in the velocity discretization). This fastens the algorithm in practice without any loss of accuracy, especially for 
uncertain problems.
\end{remark}

Now we apply an IMEX-RK scheme to system (\ref{15}) where $(f_1,g_1)^T$ is evaluated explicitly and $(f_2, g_2)^T$ implicitly, then we obtain
\begin{subequations}\label{20}
\begin{align}
&\hat{\bold r}^{n+1}=\hat{\bold r}^n+\Delta t\sum_{k=1}^s\tilde{b}_kf_1(\hat{\bold R}^k,\hat{\bold J}^k)+\Delta t\sum_{k=1}^sb_kf_2(\hat{\bold R}^k),\\
&\hat{\bold j}^{n+1}=\hat{\bold j}^n+\Delta t\sum_{k=1}^s\tilde{b}_kg_1(\hat{\bold R}^k)+\Delta t\sum_{k=1}^sb_kg_2(\hat{\bold R}^k,\hat{\bold J}^k),
\end{align}
\end{subequations}
where the internal stages are
\begin{subequations}\label{21}
\begin{align}
&\hat{\bold R}^k=\hat{\bold r}^n+\Delta t\sum_{l=1}^{k-1}\tilde{a}_{kl}f_1(\hat{\bold R}^l,\hat{\bold J}^l)+\Delta t\sum_{l=1}^ka_{kl}f_2(\hat{\bold R}^l),\\
&\hat{\bold J}^k=\hat{\bold j}^n+\Delta t\sum_{l=1}^{k-1}\tilde{a}_{kl}g_1(\hat{\bold R}^l)+\Delta t\sum_{l=1}^ka_{kl}g_2(\hat{\bold R}^l,\hat{\bold J}^l).
\end{align}
\end{subequations}

It is obvious that the scheme is characterized by the $s\times s$ matrices
\begin{equation}\label{22}
\tilde A=(\tilde a_{ij}),A=(a_{ij})
\end{equation}
and the vectors $\tilde b, b\in\mathbb R^s$, which can be represented by a double table tableau in the usual Butcher notation

\begin{center}
\begin{TAB}(r,0.5cm,0.5cm)[5pt]{c|c}{c|c}
$\tilde c$& $\tilde A$\\
& $\tilde b^T$\\
\end{TAB}, \ \ \ \ \ \ \ \ \
\begin{TAB}(r,0.5cm,0.5cm)[5pt]{c|c}{c|c}
$c$& $A$\\
& $b^T$\\
\end{TAB}.
\end{center}

The coefficients $\tilde c$ and $c$ depend on the explicit part of the scheme:
\begin{equation}\label{23}
\tilde c_i=\sum_{j=1}^{i-1}\tilde a_{ij}, \ \ c_i=\sum_{j=1}^ia_{ij}\,.
\end{equation}

In the literature, there are two main different types of IMEX-RK schemes characterized by the structure of the matrix $A$. We are interested in the IMEX-RK method of type $A$ (see \cite{boscarino2013implicit}) where the matrix $A$ is invertible, which is not sensitive to the initial data. Of course, other IMEX-RK may be used as well, we refer to \cite{boscarino2013implicit} for further examples.

As an example, we list the SSP$(3,3,2)$ scheme, which is the second order IMEX scheme we are going to use in the numerical results:
\begin{equation}
\begin{tabular}{c|c c c} 
0&0&0&0\\ 
1/2&1/2&0&0\\
1&1/2&1/2&0\\
\hline
 &1/3&1/3&1/3
\end{tabular},\qquad
\begin{tabular}{c|c c c} 
1/4&1/4&0&0\\ 
1/4&0&1/4&0\\
1&1/3&1/3&1/3\\
\hline
 &1/3&1/3&1/3
\end{tabular}.
\label{eq:SSP2}
\end{equation}

To obtain $\hat{\bold R}^k$ in each internal stage of (\ref{21}), one needs the quantity $\hat{\boldsymbol\Rho}^k=\langle\hat{\bold R}^k\rangle=\int_0^1\hat{\bold R}^kdv$ in the implicit part $f_2(\hat{\bold R}_k)$. This quantity can be implemented by the following procedure.  Suppose one has computed $\hat{\bold R}^l$ for $l=1,\cdots, k-1$, then according to (\ref{21}a):
\begin{equation}\label{24}
\begin{aligned}
\hat{\bold R}^k&=\hat{\bold r}^n+\Delta t\sum_{l=1}^{k-1}\left[\tilde{a}_{kl}f_1(\hat{\bold R}^l,\hat{\bold J}^l)+a_{kl}f_2(\hat{\bold R}^l)\right]+\Delta t\, a_{kk}\left(-\frac{1}{\varepsilon^2}\bold S(x)(\hat{\bold R}^k-\hat{\boldsymbol \Rho}^k)+\frac{\mu}{3}\partial_x(\bold S^{-1}\partial_x\hat{\boldsymbol \Rho}^k)\right)\\
&=\overline{\hat {\bold R}}^{k-1}+\Delta t\, a_{kk}\left(-\frac{1}{\varepsilon^2}\bold S(x)(\hat{\bold R}^k-\hat{\boldsymbol \Rho}^k)+\frac{\mu}{3}\partial_x(\bold S^{-1}\partial_x\hat{\boldsymbol \Rho}^k)\right).\\\end{aligned}
\end{equation}
Here $\overline{\hat {\bold R}}^{k-1}$ represents all contributions in (\ref{24}) from the first $k-1$ stages. Now one takes $\langle\cdot\rangle$ on both sides of (\ref{24}) so that $(\hat{\bold R}^k-\hat{\boldsymbol \Rho}^k)$ is cancelled out on the right hand side and one can obtain $\hat{\boldsymbol \Rho}^k$ from the following diffusion equation in implicit form:
\begin{equation}\label{25}
\hat{\boldsymbol \Rho}^k-\Delta t a_{kk}\frac{\mu}{3}\partial_x(\bold S^{-1}\partial_x\hat{\boldsymbol \Rho}^k)=\langle \overline{\hat {\bold R}}^{k-1}\rangle.
\end{equation}
 Then it is plugged  back to (\ref{21}a) in order to compute $\hat{\bold R}^k$.
 

\subsection{The Space Discretization}
To obtain the second order accuracy, we apply the upwind TVD scheme in the explicit transport part and center difference for other second derivatives. During each internal stage (\ref{21}),
\begin{subequations}\label{26}
\begin{align}
\hat{\bold R}^k_i=&\hat{\bold r}^n_i+\Delta t\sum_{l=1}^{k-1}\tilde a_{kl}\left\{-\frac{v}{2\Delta x}(\hat{\bold J}^l_{i+1}-\hat{\bold J}^l_{i-1})+\frac{v\phi^{1/2}}{2\Delta x}(\hat{\bold R}^l_{i+1}-2\hat{\bold R}^l_i+\hat{\bold R}^l_{i-1})\right.\nonumber\\
&-\frac{v\phi^{1/2}}{4}(\boldsymbol\gamma^l_i-\boldsymbol\gamma^l_{i-1}+\boldsymbol\beta^l_{i+1}-\boldsymbol\beta^l_i)\nonumber\\
&\left.-\frac{\mu}{3\Delta x^2}\left[\bold S_{i+1/2}^{-1}(\hat{\boldsymbol \Rho}^l_{i+1}-\hat{\boldsymbol \Rho}^l_{i})-\bold S_{i-1/2}^{-1}(\hat{\boldsymbol \Rho}^l_{i}-\hat{\boldsymbol \Rho}^l_{i-1})\right]\right\}\nonumber\\
&+\Delta t \sum_{l=1}^ka_{kl}\left\{-\frac{1}{\varepsilon^2}\bold S_i(\hat{\bold R}^l_i-\hat{\boldsymbol \Rho}^l_i)\right.\nonumber\\
&\left.+\frac{\mu}{3\Delta x^2}\left[\bold S_{i+1/2}^{-1}(\hat{\boldsymbol \Rho}^l_{i+1}-\hat{\boldsymbol \Rho}^l_{i})-\bold S_{i-1/2}^{-1}(\hat{\boldsymbol \Rho}^l_{i}-\hat{\boldsymbol \Rho}^l_{i-1})\right]\right\},\\
\hat{\bold J}^k_i=&\hat{\bold j}^n_i+\Delta t\sum_{l=1}^{k-1}\tilde a_{kl}\left\{-\frac{v\phi}{2\Delta x}(\hat{\bold R}^l_{i+1}-\hat{\bold R}^l_{i-1})+\frac{v\phi^{1/2}}{2\Delta x}(\hat{\bold J}^l_{i+1}-2\hat{\bold J}^l_i+\hat{\bold J}^l_{i-1})\right.\nonumber\\
&\left.-\frac{v\phi}{4}(\boldsymbol\gamma^l_i-\boldsymbol\gamma^l_{i-1}-\boldsymbol\beta^l_{i+1}+\boldsymbol\beta^l_i)\right\}\nonumber\\
&-\Delta t\sum_{l=1}^ka_{kl}\frac{1}{\varepsilon^2}\left\{\bold S_i\hat{\bold J}^l_i+(1-\varepsilon^2\phi)v\frac{\hat{\bold R}_{i+1}^l-\hat{\bold R}_{i-1}^l}{2\Delta x}\right\},
\end{align}
\end{subequations}
where
\begin{subequations}\label{27}
\begin{align}
\boldsymbol \gamma^l_i=&\frac{1}{\Delta x}\text{minmod}\left(\hat{\bold R}^l_{i+1}+\phi^{-1/2}\hat{\bold J}^l_{i+1}-\hat{\bold R}^l_i-\phi^{-1/2}\hat{\bold J}^l_i,\right. \nonumber\\
&\left.\hat{\bold R}^l_{i}+\phi^{-1/2}\hat{\bold J}^l_{i}-\hat{\bold R}^l_{i-1}-\phi^{-1/2}\hat{\bold J}^l_{i-1}\right),\\
\boldsymbol \beta^l_i=&\frac{1}{\Delta x}\text{minmod}\left(\hat{\bold R}^l_{i+1}-\phi^{-1/2}\hat{\bold J}^l_{i+1}-\hat{\bold R}^l_i+\phi^{-1/2}\hat{\bold J}^l_i,\right.\nonumber\\
&\left.\hat{\bold R}^l_{i}-\phi^{-1/2}\hat{\bold J}^l_{i}-\hat{\bold R}^l_{i-1}+\phi^{-1/2}\hat{\bold J}^l_{i-1}\right).
\end{align}
\end{subequations}

Since $\hat{\boldsymbol \Rho}^k$ can be obtained by (\ref{25}), one can get fully discretized $\hat{\boldsymbol \Rho}^k_i$ as following:
\begin{equation}\label{28}
\begin{aligned}
\hat{\boldsymbol \Rho}^k_i-\Delta ta_{kk}\frac{\mu}{3\Delta x^2}\left[\bold S_{i+1/2}^{-1}(\hat{\boldsymbol \Rho}^k_{i+1}-\hat{\boldsymbol \Rho}^k_{i})-\bold S_{i-1/2}^{-1}(\hat{\boldsymbol \Rho}^k_{i}-\hat{\boldsymbol \Rho}^k_{i-1})\right]=&\langle\overline{\hat {\bold R}}_i^{k-1}\rangle.
\end{aligned}
\end{equation}

By inverting a tridiagonal matrix, one can obtain $\hat{\boldsymbol \Rho}^k_i$. Then using (\ref{28}), the fully discretized $\hat{\bold R}^k_i$ is obtained and thus $\hat{\bold J}^k_i$ subsequently as following:
\begin{subequations}\label{29}
\begin{align}
\left(\bold I+\frac{a_{kk}\Delta t}{\varepsilon^2}\bold S_i\right)\hat{\bold R}^k_i=&\hat{\bold r}^n_i+\Delta t\sum_{l=1}^{k-1}\tilde a_{kl}\left\{-\frac{v}{2\Delta x}(\hat{\bold J}^l_{i+1}-\hat{\bold J}^l_{i-1})+\frac{v\phi^{1/2}}{2\Delta x}(\hat{\bold R}^l_{i+1}-2\hat{\bold R}^l_i+\hat{\bold R}^l_{i-1})\right.\nonumber\\
&-\frac{v\phi^{1/2}}{4}(\boldsymbol\gamma^l_i-\boldsymbol\gamma^l_{i-1}+\boldsymbol\beta^l_{i+1}-\boldsymbol\beta^l_i)\nonumber\\
&\left.-\frac{\mu}{3\Delta x^2}\left[\bold S_{i+1/2}^{-1}(\hat{\boldsymbol \Rho}^l_{i+1}-\hat{\boldsymbol \Rho}^l_{i})-\bold S_{i-1/2}^{-1}(\hat{\boldsymbol \Rho}^l_{i}-\hat{\boldsymbol \Rho}^l_{i-1})\right]\right\}\nonumber\\
&+\Delta t \sum_{l=1}^{k-1}a_{kl}\left\{-\frac{1}{\varepsilon^2}\bold S_i(\hat{\bold R}^l_i-\hat{\boldsymbol \Rho}^l_i)\right.\nonumber\\
&\left.+\frac{\mu}{3\Delta x^2}\left[\bold S_{i+1/2}^{-1}(\hat{\boldsymbol \Rho}^l_{i+1}-\hat{\boldsymbol \Rho}^l_{i})-\bold S_{i-1/2}^{-1}(\hat{\boldsymbol \Rho}^l_{i}-\hat{\boldsymbol \Rho}^l_{i-1})\right]\right\}\nonumber\\
&+\Delta t a_{kk}\left\{\frac{1}{\varepsilon^2}\bold S_i\hat{\boldsymbol \Rho}^k_i+\frac{\mu}{3\Delta x^2}\left[\bold S_{i+1/2}^{-1}(\hat{\boldsymbol \Rho}^k_{i+1}-\hat{\boldsymbol \Rho}^k_{i})\right.\right.\nonumber\\
&\left.\left.-\bold S_{i-1/2}^{-1}(\hat{\boldsymbol \Rho}^k_{i}-\hat{\boldsymbol \Rho}^k_{i-1})\right]\right\},\\
\left(1+\frac{a_{kk}\Delta t}{\varepsilon^2}\bold S_i\right)\hat{\bold J}^k_i=&\hat{\bold j}^n_i+\Delta t\sum_{l=1}^{k-1}\tilde a_{kl}\left\{-\frac{v\phi}{2\Delta x}(\hat{\bold R}^l_{i+1}-\hat{\bold R}^l_{i-1})+\frac{v\phi^{1/2}}{2\Delta x}(\hat{\bold J}^l_{i+1}-2\hat{\bold J}^l_i+\hat{\bold J}^l_{i-1})\right.\nonumber\\
&\left.-\frac{v\phi}{4}(\boldsymbol\gamma^l_i-\boldsymbol\gamma^l_{i-1}+\boldsymbol\beta^l_{i+1}-\boldsymbol\beta^l_i)\right\}\nonumber\\
&-\Delta t\sum_{l=1}^{k-1}a_{kl}\frac{1}{\varepsilon^2}\left\{\bold S_i\hat{\bold J}^l_i+(1-\varepsilon^2\phi)v\frac{\hat{\bold R}_{i+1}^l-\hat{\bold R}_{i-1}^l}{2\Delta x}\right\}\nonumber\\
&-\Delta ta_{kk}\left(\frac{1-\varepsilon^2\phi}{\varepsilon^2}\right)v\frac{\hat{\bold R}_{i+1}^k-\hat{\bold R}_{i-1}^k}{2\Delta x}.
\end{align}
\end{subequations}

After calculating all $\hat{\bold R}^k_i$ and $\hat{\bold J}^k_i$ for $k=1,\cdots, s$, one can update $\hat{\bold r}^{n+1}_i$ and $\hat{\bold j}^{n+1}_i$ in (\ref{20}),
\begin{subequations}\label{30}
\begin{align}
\hat{\bold r}^{n+1}_i=&\hat{\bold r}^n_i+\Delta t\sum_{k=1}^{s}\tilde b_{k}\left\{-\frac{v}{2\Delta x}(\hat{\bold J}^k_{i+1}-\hat{\bold J}^k_{i-1})+\frac{v\phi^{1/2}}{2\Delta x}(\hat{\bold R}^k_{i+1}-2\hat{\bold R}^k_i+\hat{\bold R}^k_{i-1})\right.\nonumber\\
&-\frac{v\phi^{1/2}}{4}(\boldsymbol\gamma^k_i-\boldsymbol\gamma^k_{i-1}+\boldsymbol\beta^k_{i+1}-\boldsymbol\beta^k_i)\nonumber\\
&\left.-\frac{\mu}{3\Delta x^2}\left[\bold S_{i+1/2}^{-1}(\hat{\boldsymbol \Rho}^k_{i+1}-\hat{\boldsymbol \Rho}^k_{i})-\bold S_{i-1/2}^{-1}(\hat{\boldsymbol \Rho}^k_{i}-\hat{\boldsymbol \Rho}^k_{i-1})\right]\right\}\nonumber\\
&+\Delta t \sum_{k=1}^sb_{k}\left\{-\frac{1}{\varepsilon^2}\bold S_i(\hat{\bold R}^k_i-\hat{\boldsymbol \Rho}^k_i)\right.\nonumber\\
&\left.+\frac{\mu}{3\Delta x^2}\left[\bold S_{i+1/2}^{-1}(\hat{\boldsymbol \Rho}^k_{i+1}-\hat{\boldsymbol \Rho}^k_{i})-\bold S_{i-1/2}^{-1}(\hat{\boldsymbol \Rho}^k_{i}-\hat{\boldsymbol \Rho}^k_{i-1})\right]\right\},\\
\hat{\bold j}^{n+1}_i=&\hat{\bold j}^n_i+\Delta t\sum_{k=1}^{s}\tilde b_{k}\left\{-\frac{v\phi}{2\Delta x}(\hat{\bold R}^k_{i+1}-\hat{\bold R}^k_{i-1})+\frac{v\phi^{1/2}}{2\Delta x}(\hat{\bold J}^k_{i+1}-2\hat{\bold J}^k_i+\hat{\bold J}^k_{i-1})\right.\nonumber\\
&\left.-\frac{v\phi}{4}(\boldsymbol\gamma^k_i-\boldsymbol\gamma^k_{i-1}+\boldsymbol\beta^k_{i+1}-\boldsymbol\beta^k_i)\right\}\nonumber\\
&-\Delta t\sum_{k=1}^kb_{k}\frac{1}{\varepsilon^2}\left\{\bold S_i\hat{\bold J}^k_i+(1-\varepsilon^2\phi)v\frac{\hat{\bold R}_{i+1}^k-\hat{\bold R}_{i-1}^k}{2\Delta x}\right\},
\end{align}
\end{subequations}
where $\boldsymbol \gamma^k_i$ and $\boldsymbol \beta^k_i$ are defined the same as in (\ref{27}).

We choose \cite{boscarino2013implicit}
\begin{equation}\label{31}
\mu=\exp(-\varepsilon^2/\Delta x).
\end{equation}
Thus, for large value of $\varepsilon$, (e.g., $\varepsilon=1$), one can
avoid the loss of accuracy caused by adding and subtracting the penalty term; for very small value of $\varepsilon$, (e.g., $\varepsilon\to 0$), $\mu\to1$.

Before continuing some observations are appropriate.

\begin{remark}~
\begin{itemize}
\item
The velocity space is discretized by the so-called discrete ordinate method
\cite{CZ}  which corresponds to use the Gauss-Legendre quadrature for the velocity integrals.
\item
Notice that the matrix $\bold S(x)$ can be precomputed. For each stage of (\ref{20}), a large tridiagonal matrix has to be inverted. Since it is sparse and can be precomputed, fast algorithms can be applied. 
\item
To get the boundary conditions for $\hat{\bold r}$ and $\hat{\bold j}$ we refer to \cite{jin2000uniformly} for details.
\end{itemize}
\end{remark}

\subsection{The sAP property}
Denote
$$\begin{aligned}
f_1(\hat{\bold R}^l_i,\hat{\bold J}^l_i)=&-\frac{v}{2\Delta x}(\hat{\bold J}^l_{i+1}-\hat{\bold J}^l_{i-1})+\frac{v\phi^{1/2}}{2\Delta x}(\hat{\bold R}^l_{i+1}-2\hat{\bold R}^l_i+\hat{\bold R}^l_{i-1})\\
&-\frac{v\phi^{1/2}}{4}(\boldsymbol\gamma^l_i-\boldsymbol\gamma^l_{i-1}+\boldsymbol\beta^l_{i+1}-\boldsymbol\beta^l_i)\\
&-\frac{\mu}{3\Delta x^2}\left[\bold S_{i+1/2}^{-1}(\hat{\boldsymbol \Rho}^l_{i+1}-\hat{\boldsymbol \Rho}^l_{i})-\bold S_{i-1/2}^{-1}(\hat{\boldsymbol \Rho}^l_{i}-\hat{\boldsymbol \Rho}^l_{i-1})\right],\\
f_2(\hat{\bold R}^l_i)=&-\frac{1}{\varepsilon^2}\bold S_i(\hat{\bold R}^l_i-\hat{\boldsymbol \Rho}^l_i)\nonumber\\
&+\frac{\mu}{3\Delta x^2}\left[\bold S_{i+1/2}^{-1}(\hat{\boldsymbol \Rho}^l_{i+1}-\hat{\boldsymbol \Rho}^l_{i})-\bold S_{i-1/2}^{-1}(\hat{\boldsymbol \Rho}^l_{i}-\hat{\boldsymbol \Rho}^l_{i-1})\right],\\
g_1(\hat{\bold R}^l_i)=&-\frac{v\phi}{2\Delta x}(\hat{\bold R}^l_{i+1}-\hat{\bold R}^l_{i-1})+\frac{v\phi^{1/2}}{2\Delta x}(\hat{\bold J}^l_{i+1}-2\hat{\bold J}^l_i+\hat{\bold J}^l_{i-1})\\
&-\frac{v\phi}{4}(\boldsymbol\gamma^l_i-\boldsymbol\gamma^l_{i-1}-\boldsymbol\beta^l_{i+1}+\boldsymbol\beta^l_i),\\
g_2(\hat{\bold R}^l_i,\hat{\bold J}^l_i)=&\frac{1}{\varepsilon^2}\left[\bold S_i\hat{\bold J}^l_i+(1-\varepsilon^2\phi)v\frac{\hat{\bold R}_{i+1}^l-\hat{\bold R}_{i-1}^l}{2\Delta x}\right].
\end{aligned}$$

From (\ref{26}) one gets
\begin{subequations}\label{32}
\begin{align}
\begin{pmatrix}
\hat{\bold R}_i^1\\
\hat{\bold R}_i^2\\
\vdots\\
\hat{\bold R}_i^s
\end{pmatrix}=&
\begin{pmatrix}
\hat{\bold r}_i^n\\
\hat{\bold r}_i^n\\
\vdots\\
\hat{\bold r}_i^n
\end{pmatrix}+\Delta t\begin{pmatrix}
0\\
\tilde a_{21}f_1(\hat{\bold R}_i^1,\hat{\bold J}_i^1)\\
\vdots\\
\sum_{l=1}^{s-1}\tilde a_{sl}f_1(\hat{\bold R}_i^{l},\hat{\bold J}_i^{l})
\end{pmatrix}+\Delta t\bold A
\begin{pmatrix}
f_2(\hat{\bold R}_i^1)\\
f_2(\hat{\bold R}_i^2)\\
\vdots\\
f_2(\hat{\bold R}_i^{s})
\end{pmatrix},\\
\begin{pmatrix}
\hat{\bold J}_i^1\\
\hat{\bold J}_i^2\\
\vdots\\
\hat{\bold J}_i^s
\end{pmatrix}=&
\begin{pmatrix}
\hat{\bold j}_i^n\\
\hat{\bold j}_i^n\\
\vdots\\
\hat{\bold j}_i^n
\end{pmatrix}+\Delta t\begin{pmatrix}
0\\
\tilde a_{21}g_1(\hat{\bold R}_i^1)\\
\vdots\\
\sum_{l=1}^{s-1}\tilde a_{sl}g_1(\hat{\bold R}_i^{l})
\end{pmatrix}+\Delta t\bold A\begin{pmatrix}
g_2(\hat{\bold R}_i^1,\hat{\bold J}_i^1)\\
g_2(\hat{\bold R}_i^2,\hat{\bold J}_i^2)\\
\vdots\\
g_2(\hat{\bold R}_i^{s},\hat{\bold J}_i^s)
\end{pmatrix},
\end{align}
\end{subequations}
where 
\begin{equation}\label{33}
\bold A_{K(i-1)+1:Ki,K(j-1)+1:Kj}=A_{i,j}\bold I_{K\times K}, \ \ \ \bold I_{K\times K} \ \ \text{is the } \ K\times K \ \text{identity matrix},
\end{equation}
and $A$ is defined in (\ref{21}). Denote $\bold W$ as the inverse matrix of $\bold A$, then one obtains from (\ref{32})
\begin{subequations}\label{34}
\begin{align}
&\Delta t
\begin{pmatrix}
f_2(\hat{\bold R}_i^1)\\
f_2(\hat{\bold R}_i^2)\\
\vdots\\
f_2(\hat{\bold R}_i^{s})
\end{pmatrix}=\bold W\left[\begin{pmatrix}
\hat{\bold R}_i^1\\
\hat{\bold R}_i^2\\
\vdots\\
\hat{\bold R}_i^s
\end{pmatrix}-
\begin{pmatrix}
\hat{\bold r}_i^n\\
\hat{\bold r}_i^n\\
\vdots\\
\hat{\bold r}_i^n
\end{pmatrix}-\Delta t\begin{pmatrix}
0\\
\tilde a_{21}f_1(\hat{\bold R}_i^1,\hat{\bold J}_i^1)\\
\vdots\\
\sum_{l=1}^{s-1}\tilde a_{sl}f_1(\hat{\bold R}_i^{l},\hat{\bold J}_i^{l})
\end{pmatrix}\right],\\
&\Delta t\begin{pmatrix}
g_2(\hat{\bold R}_i^1,\hat{\bold J}_i^1)\\
g_2(\hat{\bold R}_i^2,\hat{\bold J}_i^2)\\
\vdots\\
g_2(\hat{\bold R}_i^{s},\hat{\bold J}_i^s)
\end{pmatrix}=\bold W\left[\begin{pmatrix}
\hat{\bold J}_i^1\\
\hat{\bold J}_i^2\\
\vdots\\
\hat{\bold J}_i^s
\end{pmatrix}-
\begin{pmatrix}
\hat{\bold j}_i^n\\
\hat{\bold j}_i^n\\
\vdots\\
\hat{\bold j}_i^n
\end{pmatrix}-\Delta t\begin{pmatrix}
0\\
\tilde a_{21}g_1(\hat{\bold R}_i^1)\\
\vdots\\
\sum_{l=1}^{s-1}\tilde a_{sl}g_1(\hat{\bold R}_i^{l})
\end{pmatrix}\right],
\end{align}
\end{subequations}

Since $\bold W$ has the same structure as $\bold A$, $\bold W$ should be a lower  triangular matrix with entries
\begin{equation}\label{35}
\bold W_{K(i-1)+1:Ki,K(j-1)+1:Kj}=\omega_{i,j}\bold I_{K\times K},
\end{equation}
where $W=(\omega_{i,j})$ is the inverse of the lower triangular matrix $A$ in (\ref{21}).

Then one can rewrite (\ref{34}) as 
\begin{subequations}\label{36}
\begin{align}
&\Delta tf_2(\hat{\bold R}_i^k)=\sum_{l=1}^k\omega_{kl}\left[\hat{\bold R}_i^l-\hat{\bold r}_i^n-\Delta t\sum_{l=1}^{k-1}\tilde{a}_{kl}f_1(\hat{\bold R}_i^l,\hat{\bold J}_i^l)\right],\\
&\Delta tg_2(\hat{\bold R}_i^k,\hat{\bold J}_i^k)=\sum_{l=1}^k\omega_{kl}\left[\hat{\bold J}_i^l-\hat{\bold j}_i^n-\Delta t\sum_{l=1}^{k-1}\tilde a_{kl}g_1(\hat{\bold R}_i^l)\right],
\end{align}
\end{subequations}

More explicitly,
\begin{subequations}\label{37}
\begin{align}
&\Delta t\left\{-\frac{1}{\varepsilon^2}\bold S_i(\hat{\bold R}^k_i-\hat{\boldsymbol \Rho}^k_i)+\frac{\mu}{3\Delta x^2}\left[\bold S_{i+1/2}^{-1}(\hat{\boldsymbol \Rho}^k_{i+1}-\hat{\boldsymbol \Rho}^k_{i})-\bold S_{i-1/2}^{-1}(\hat{\boldsymbol \Rho}^k_{i}-\hat{\boldsymbol \Rho}^k_{i-1})\right]\right\}\nonumber\\
=&\sum_{l=1}^k\omega_{kl}\left\{\hat{\bold R}_i^l-\hat{\bold r}_i^n-\Delta t\sum_{l=1}^{k-1}\tilde{a}_{kl}\left[-\frac{v}{2\Delta x}(\hat{\bold J}^l_{i+1}-\hat{\bold J}^l_{i-1})+\frac{v\phi^{1/2}}{2\Delta x}(\hat{\bold R}^l_{i+1}-2\hat{\bold R}^l_i+\hat{\bold R}^l_{i-1})\right.\right.\nonumber\\
&\left.\left.-\frac{v\phi^{1/2}}{4}(\boldsymbol\gamma^l_i-\boldsymbol\gamma^l_{i-1}+\boldsymbol\beta^l_{i+1}-\boldsymbol\beta^l_i)-\frac{\mu}{3\Delta x^2}\left(\bold S_{i+1/2}^{-1}(\hat{\boldsymbol \Rho}^l_{i+1}-\hat{\boldsymbol \Rho}^l_{i})-\bold S_{i-1/2}^{-1}(\hat{\boldsymbol \Rho}^l_{i}-\hat{\boldsymbol \Rho}^l_{i-1})\right)\right]\right\},\\
&\Delta t\left\{\frac{1}{\varepsilon^2}\left[\bold S_i\hat{\bold J}^k_i+(1-\varepsilon^2\phi)v\frac{\hat{\bold R}_{i+1}^k-\hat{\bold R}_{i-1}^k}{2\Delta x}\right]\right\}\nonumber\\
=&\sum_{l=1}^k\omega_{kl}\left\{\hat{\bold J}_i^l-\hat{\bold j}_i^n-\Delta t \sum_{l=1}^{k-1}\tilde a_{kl}\left[-\frac{v\phi}{2\Delta x}(\hat{\bold R}^l_{i+1}-\hat{\bold R}^l_{i-1})+\frac{v\phi^{1/2}}{2\Delta x}(\hat{\bold J}^l_{i+1}-2\hat{\bold J}^l_i+\hat{\bold J}^l_{i-1})\right.\right.\nonumber\\
&\left.\left.-\frac{v\phi}{4}(\boldsymbol\gamma^l_i-\boldsymbol\gamma^l_{i-1}-\boldsymbol\beta^l_{i+1}+\boldsymbol\beta^l_i)\right]\right\}.
\end{align}
\end{subequations}

Thus, setting $\varepsilon\to 0$, since $\bold S_i$ is non-singular, one obtains
\begin{subequations}\label{38}
\begin{align}
&\hat{\bold R}_i^k=\hat{\boldsymbol \Rho}_i^k,\\
&\hat{\bold J}_i^k=-v\bold S_i^{-1}\frac{\hat{\bold R}_{i+1}^k-\hat{\bold R}_{i-1}^k}{2\Delta x}.
\end{align}
\end{subequations}

Inserting this back to (\ref{30}a) and letting $\varepsilon\to 0$,
\begin{equation}\label{39}
\hat{\bold r}_i^{n+1}=\hat{\bold r}_i^n+\Delta t\sum_{k=1}^s\tilde{b}_k\hat f_1(\hat{\bold R}_i^k)+\Delta t\sum_{k=1}^sb_k\hat f_2(\hat{\bold R}_i^k),
\end{equation}
where
\begin{subequations}\label{40}
\begin{align}
\hat f_1(\hat{\bold R}_i^k)=&v^2\frac{1}{4\Delta x^2}\left[\bold S_{i+1}^{-1}(\hat{\bold R}_{i+2}^k-\hat{\bold R}_{i}^k)-\bold S_{i-1}^{-1}(\hat{\bold R}_{i}^k-\hat{\bold R}_{i-2}^k)\right]\\
&-\frac{1}{3\Delta x^2}\left[\bold S_{i+1/2}^{-1}(\hat{\bold R}_{i+1}^k-\hat{\bold R}_i^k)-\bold S_{i-1/2}^{-1}(\hat{\bold R}_i^k-\hat{\bold R}_{i-1}^k)\right],\\
\hat f_2(\hat{\bold R}_i^k)=&\frac{1}{3\Delta x^2}\left[\bold S_{i+1/2}^{-1}(\hat{\bold R}_{i+1}^k-\hat{\bold R}_i^k)-\bold S_{i-1/2}^{-1}(\hat{\bold R}_i^k-\hat{\bold R}_{i-1}^k)\right].
\end{align}
\end{subequations}
Integrating above over $v$, one gets
\begin{equation}\label{41}
\hat{\boldsymbol \rho}_i^{n+1}=\hat{\boldsymbol \rho}_i^n+\frac{\Delta t}{3\Delta x^2}\sum_{k=1}^sb_k\left[\bold S_{i+1/2}^{-1}(\hat{\boldsymbol \Rho}_{i+1}^k-\hat{\boldsymbol \Rho}_i^k)-\bold S_{i-1/2}^{-1}(\hat{\boldsymbol \Rho}_i^k-\hat{\boldsymbol \Rho}_{i-1}^k)\right]+O(\Delta x^2),
\end{equation}
where
\begin{equation}\label{42}
\hat{\boldsymbol \Rho}_i^k=\hat{\boldsymbol \rho}_i^n+\frac{\Delta t}{3\Delta x^2}\sum_{l=1}^ka_{kl}\left[\bold S_{i+1/2}^{-1}(\hat{\boldsymbol \Rho}_{i+1}^l-\hat{\boldsymbol \Rho}_i^l)-\bold S_{i-1/2}^{-1}(\hat{\boldsymbol \Rho}_i^l-\hat{\boldsymbol \Rho}_{i-1}^l)\right],
\end{equation}
which is an {\it implicit} RK scheme (requiring inverting a tridiagonal matrix at each stage) for the limiting diffusion equation (\ref{12}) and also corresponds to (\ref{28}). Thus, the sAP property \cite{jin2015asymptotic} of the efficient IMEX R-K scheme is shown.

\section{The Radiative Heat Transfer Equation with Random Inputs}
Let $x\in D\subset \mathbb{R}^3$ be the space variable, $\Omega\in S^2$ be the direction variable, $S^2$ the unit sphere of $\mathbb{R}^3$, $z\in\mathbb{R}^d  (d\geq1)$ be the random variable and $t\in\mathbb{R}^+$ the time.

We denote by $I=I(x,\Omega,z,t)$ the radiative intensity  and by $\theta(x,z,t)$ the material temperature. Introducing the Knudsen number $\varepsilon$, the radiative heat transfer equations in nondimensional form are
\begin{subequations} \label{43}
\begin{align}
&\varepsilon^2M\partial_t I+\varepsilon\Omega\cdot\nabla_xI=B(\theta)-I\\
&\varepsilon^2\partial_t \theta=\varepsilon^2\Delta_x\theta-(B(\theta)-\langle I\rangle)
\end{align}
\end{subequations}
with the total intensity
\begin{equation} \label{44}
\langle I \rangle(x,z,t)=\frac{1}{|S|^2}\int_{S^2}I(x,\Omega,z,t)d\Omega,
\end{equation}
and the black body intensity
\begin{equation} \label{45}
B(\theta)=\sigma \theta^4,
\end{equation}
where $M$ is the Mach number ($=1$ for this paper) and $\sigma=\sigma(x, z)>0$ is the cross-section depending on the space variable and the random variable.

The initial conditions and reflection transmission boundary conditions are prescribed as following:
\begin{equation} \label{46}
\begin{split}
&I.C. \ \left\{
\begin{aligned}
&I(x,\Omega,z,0)=I_I(x,\Omega,z),\\
&\theta(x,z,0)=\theta_I(x,z)
\end{aligned}
\right.
\\
&B.C. \ \left\{
\begin{aligned}
&I(\hat{x},\Omega,z,t)=\alpha(n(\hat{x})\cdot \Omega)I(\hat{x},\Omega',z,t)+[1-\alpha(n(\hat{x})\cdot \Omega)]I_B(\hat{x},\Omega,z,t),\\
&\ \ \ \ \ \ \ \ \ \ \ \ \ \ \ \ \ \ \text{for} \ n(\hat{x})\cdot \Omega<0\\
&\theta(\hat{x},z,t)=\theta_B(\hat{x},z,t)
\end{aligned}
\right.
\end{split}
\end{equation}
where $\hat{x}\in\partial D$ with outward unit normal $n(\hat{x})$ and $\Omega'=\Omega-2n(\hat{x})(n(\hat{x})\cdot\Omega)$ is the reflection of $\Omega$ on the tangent plane to $\partial D$. The reflectivity $\alpha,\ 0\leq\alpha\leq1$, depends on the incidence angle.

\subsection{An Even-Odd Decomposition}
For simplicity, we consider the one-dimensional case $x\in[0,1]$ and define $v=\cos(\Omega\cdot x), v\in[-1,1]$. Thus, the angular averaging is defined as:
$$\langle f\rangle=\frac{1}{2}\int_{-1}^1f(v)dv.$$
The one-dimensional radiative heat transfer equation is
\begin{subequations} \label{47}
\begin{align}
&\varepsilon^2\partial_t I+\varepsilon v\partial_xI=B(\theta)-I\\
&\varepsilon^2\partial_t \theta=\varepsilon^2\partial_{xx}\theta-(B(\theta)-\langle I\rangle)
\end{align}
\end{subequations}
For $v>0$,
\begin{equation}\label{48}
\begin{aligned}
&\varepsilon^2\partial_t I(v)+\varepsilon v\partial_xI(v)=B(\theta)-I(v)\\
&\varepsilon^2\partial_t I(-v)-\varepsilon v\partial_xI(-v)=B(\theta)-I(-v)\\
&\varepsilon^2\partial_t \theta=\varepsilon^2\partial_{xx}\theta-(B(\theta)-\langle I\rangle)
\end{aligned}
\end{equation}

Now denote the even and odd parities
\begin{equation}\label{49}
\begin{aligned}
&r(t,x,v,z)=\frac{1}{2}(I(t,x,v,z)+I(t,x,-v,z)),\\
&j(t,x,v,z)=\frac{1}{2\varepsilon}(I(t,x,v,z)-I(t,x,-v,z)).
\end{aligned}
\end{equation}

Then (\ref{48}) becomes
\begin{subequations}\label{50}
\begin{align}
&\partial_t r+v\partial_x j=\frac{1}{\varepsilon^2}(B(\theta)-r),\\
&\partial_t j+\frac{1}{\varepsilon^2}v\partial_x r=-\frac{1}{\varepsilon^2}j,\\
&\partial_t \theta=\partial_{xx}\theta-\frac{1}{\varepsilon^2}(B(\theta)-\langle r\rangle).
\end{align}
\end{subequations}

As $\varepsilon\to 0$, (\ref{50}) gives
\begin{equation}\label{51}
B(\theta)=r=\langle r\rangle, \ \ \ j=-v\partial_x r.
\end{equation}
Applying this in (\ref{50}a), integrating over $v$ and adding to (\ref{50}c)
\begin{equation}\label{52}
\partial_t(\theta+B(\theta))=\partial_x\left[\partial_x \theta+\frac{1}{3}\partial_xB(\theta)\right],
\end{equation}
which is the same limiting diffusion equation one can derive from the Hilbert expansion as \cite{klar2001numerical}.

\subsection{The gPC-SG Formulation}
Now we deal with the uncertainty. Using gPC approximation for $r(t,x,v,z), j(t,x,v,z)$ and $\theta(t,x,z)$ and truncate at the $N$-th order,
\begin{subequations}\label{53}
\begin{align}
&r(t,x,v,z)\approx r_N(t,x,v,z)=\sum_{k=1}^K\hat{r}_k(t,x,v)\Phi_k(z)=\hat{\bold r}\cdot \boldsymbol \Phi ,\\
&j(t,x,v,z)\approx j_N(t,x,v,z)=\sum_{k=1}^K\hat{j}_k(t,x,v)\Phi_k(z)=\hat{\bold j}\cdot \boldsymbol \Phi,\\
&\theta(t,x,z)\approx \theta_N(t,x,z)=\sum_{k=1}^K\hat{\theta}_k(t,x)\Phi_k(z)=\hat{\boldsymbol \theta}\cdot \boldsymbol \Phi,
\end{align}
\end{subequations}
where $\Phi_k$ is the same as defined in (\ref{8}) and  
\begin{subequations}\label{54}
\begin{align}
&\hat{\boldsymbol r}=(\hat r_1,\hat r_2,\cdots, \hat r_K)^T,\\
&\hat{\boldsymbol j}=(\hat j_1,\hat j_2,\cdots, \hat j_K)^T,\\
&\hat{\boldsymbol \theta}=(\hat \theta_1,\hat \theta_2,\cdots, \hat \theta_K)^T,\\
&\hat{\boldsymbol \Phi}=(\hat \Phi_1,\hat \Phi_2,\cdots, \hat \Phi_K)^T.
\end{align}
\end{subequations}

Then (\ref{50}) becomes
\begin{subequations}\label{55}
\begin{align}
&\partial_t \hat{\bold r}+v\partial_x \hat{\bold j}=\frac{1}{\varepsilon^2}(\bold B(\hat{\boldsymbol \theta})-\hat{\bold r}),\\
&\partial_t \hat{\bold j}+\frac{1}{\varepsilon^2}v\partial_x\hat{\bold r}=-\frac{1}{\varepsilon^2}\hat{\bold j},\\
&\partial_t \hat{\boldsymbol \theta}=\partial_{xx}\hat{\boldsymbol \theta}-\frac{1}{\varepsilon^2}(\bold B(\hat{\boldsymbol \theta})-\langle \hat{\bold r}\rangle).
\end{align}
\end{subequations}
where $\bold B(\hat{\boldsymbol \theta})=(B_{j})_{1\leq j\leq K}$ with
\begin{equation}\label{56}
B_{j}=\int_{I_z}(\hat{\boldsymbol \theta}\cdot \boldsymbol \Phi)^4\sigma(x,z)\Phi_j(z)\pi(z)dz
\end{equation}

Similarly, the limiting diffusion equation ({\ref{52}) becomes
\begin{equation}\label{57}
\partial_t[\hat{\boldsymbol \theta}+\bold B(\hat{\boldsymbol\theta})]=\partial_x\left[(\bold I+\frac{4}{3}\bold C)\partial_x \hat{\boldsymbol \theta}\right],
\end{equation}
where $\bold I$ is a $K\times K$ identity matrix and $\bold C(\hat{\boldsymbol \theta})=(c_{ij})_{1\leq i,j\leq K}$ with
\begin{equation}\label{58}
c_{ij}=\int_{I_z}(\hat{\boldsymbol \theta}\cdot\boldsymbol \Phi)^3\sigma(x,z)\Phi_i(z)\Phi_j(z)\pi(z)dz
\end{equation}

The well-posedness of this equation under some gentle condition is proved in \cite{JinLu}.

\subsection{ An Efficient IMEX-RK Scheme}
One could apply the same relaxation method as in \cite{jin2000uniformly} for the system (\ref{55}), but it will still suffer from the parabolic CFL condition $\Delta t=O(\Delta x^2)$. So similar as in section 2.2, we rewrite (\ref{55}) as following
\begin{subequations}\label{59}
\begin{align}
&\partial_t \hat{\bold r}=-v\partial_x \hat{\bold j}+\frac{1}{\varepsilon^2}(\bold B(\hat{\boldsymbol \theta})-\hat{\bold r}),\\
&\partial_t \hat{\bold j}=-\frac{1}{\varepsilon^2}(\hat{\bold j}+v\partial_x \hat{\bold r}),\\
&\partial_t \hat{\boldsymbol \theta}=\partial_{xx}\hat{\boldsymbol\theta}-\frac{1}{\varepsilon^2}(\bold B(\hat{\boldsymbol \theta})-\langle \hat{\bold r}\rangle).
\end{align}
\end{subequations}

By adding and subtracting the term $\frac{\mu}{3}\partial_{xx}\bold B(\hat{\boldsymbol \theta})$ in (\ref{59}a) and the term $\phi v\partial_x \hat{\bold r}$ in (\ref{59}b), we reformulate the problem into its equivalent form:
\begin{subequations}\label{60}
\begin{align}
\partial_t \hat{\bold r}&=\left[-v\partial_x \hat{\bold j}-\frac{\mu}{3}\partial_{xx}\bold B(\hat{\boldsymbol \theta})\right]+\left[\frac{1}{\varepsilon^2}(\bold B(\hat{\boldsymbol \theta})-\hat{\bold r})+\frac{\mu}{3}\partial_{xx}\bold B(\hat{\boldsymbol \theta})\right],\nonumber\\
&=f_1(\hat{\bold j},\hat{\boldsymbol \theta})+f_2(\hat{\bold r},\hat{\boldsymbol \theta}),\\
\partial_t \hat{\bold j}&=-v\phi\partial_x\hat{\bold r}-\frac{1}{\varepsilon^2}[\hat{\bold j}+(1-\varepsilon^2\phi)v\partial_x \hat{\bold r}]=g_1(\hat{\bold r})+g_2(\hat{\bold r},\hat{\bold j}),\\
\partial_t \hat{\boldsymbol \theta}&=\partial_{xx}\hat{\boldsymbol \theta}-\frac{1}{\varepsilon^2}(\bold B(\hat{\boldsymbol \theta})-\langle\hat{\bold r}\rangle)=h(\hat{\boldsymbol \theta}),
\end{align}
\end{subequations}
where 
\begin{subequations}\label{61}
\begin{align}
&f_1(\hat{\bold j},\hat{\boldsymbol \theta})=-v\partial_x\hat{\bold j}-\frac{\mu}{3}\partial_{xx}\bold B(\hat{\boldsymbol \theta}),\\
&f_2(\hat{\bold r},\hat{\boldsymbol \theta})=\frac{1}{\varepsilon^2}(\bold B(\hat{\boldsymbol \theta})-\hat{\bold r})+\frac{\mu}{3}\partial_{xx}\bold B(\hat{\boldsymbol \theta}),\\
&g_1(\hat{\bold r})=-v\phi\partial_x\hat{\bold r},\\
&g_2(\hat{\bold r},\hat{\bold j})=-\frac{1}{\varepsilon^2}[\hat{\bold j}+(1-\varepsilon^2\phi)v\partial_x \hat{\bold r}].
\end{align}
\end{subequations}
The constant $\mu$ and $\phi$ are chosen the same as in (\ref{17}) and (\ref{19}) respectively.

\subsubsection{The IMEX Runge-Kutta Method}
Now we apply an IMEX-RK scheme to system (\ref{60}) where $(f_1,g_1,0)^T$ is evaluated explicitly and $(f_2, g_2,h)^T$ implicitly, then we obtain
\begin{subequations}\label{62}
\begin{align}
&\hat{\bold r}^{n+1}= \hat{\bold r}^n+\Delta t\sum_{k=1}^s\tilde{b}_kf_1(\hat{\bold J}^k,\hat{\bold \Theta}^k)+\Delta t\sum_{k=1}^sb_kf_2(\hat{\bold R}^k,\hat{\bold\Theta}^k),\\
& \hat{\bold j}^{n+1}=\hat{\bold j}^n+\Delta t\sum_{k=1}^s\tilde{b}_kg_1(\hat{\bold R}^k)+\Delta t\sum_{k=1}^sb_kg_2(\hat{\bold R}^k,\hat{\bold J}^k),\\
&\hat{\boldsymbol \theta}^{n+1}=\hat{\boldsymbol \theta}^n+\Delta t\sum_{k=1}^sb_kh(\hat{\boldsymbol\Theta}^k),
\end{align}
\end{subequations}
where the internal stages are
\begin{subequations}\label{63}
\begin{align}
&\hat{\bold R}^k=\hat{\bold r}^n+\Delta t\sum_{l=1}^{k-1}\tilde{a}_{kl}f_1(\hat{\bold J}^l,\hat{\boldsymbol\Theta}^l)+\Delta t\sum_{l=1}^ka_{kl}f_2(\hat{\boldsymbol R}^l,\hat{\boldsymbol \Theta}^l),\\
&\hat{\bold J}^k=\hat{\bold j}^n+\Delta t\sum_{l=1}^{k-1}\tilde{a}_{kl}g_1(\hat{\bold R}^l)+\Delta t\sum_{l=1}^ka_{kl}g_2(\hat{\bold R}^l,\hat{\bold J}^l),\\
&\hat{\boldsymbol \Theta}^k=\hat{\boldsymbol\theta}^n+\Delta t\sum_{l=1}^ka_{kl}h(\hat{\boldsymbol\Theta}^l).
\end{align}
\end{subequations}
Notice that for each internal stage, one has to compute $\hat{\boldsymbol \Theta}_k$ first. Then $\hat{\bold R}_k$ and $\hat{\bold J}_k$ can be solved similarly as before. 
\begin{equation}\label{64}
\hat{\boldsymbol \Theta}^k=\hat{\boldsymbol \theta}^n+\Delta t\sum_{l=1}^{k}a_{kl}\left(\partial_{xx}\hat{\boldsymbol \Theta}^l-\frac{1}{\varepsilon^2}(\bold B(\hat{\boldsymbol\Theta}^l)-\langle \hat{\boldsymbol R}^l\rangle)\right).
\end{equation}

Take $\langle\cdot \rangle$ on (\ref{63}a),
\begin{equation}\label{65}
\begin{aligned}
\langle \hat{\bold R}^k\rangle=&\langle \hat{\bold r}^n +\Delta t\sum_{l=1}^{k-1}\tilde{a}_{kl}f_1(\hat{\bold J}^l,\hat{\bold \Theta}^l)+\Delta t\sum_{l=1}^{k-1}a_{kl}f_2(\hat{\bold R}^l,\hat{\bold \Theta}^l)\rangle+\Delta ta_{kk}\langle f_2(\hat{\bold R}^k,\hat{\bold \Theta}^k)\rangle\\
=&\langle \hat{\bold r}^n +\Delta t\sum_{l=1}^{k-1}\tilde{a}_{kl}f_1(\hat{\bold J}^l,\hat{\bold \Theta}^l)+\Delta t\sum_{l=1}^{k-1}a_{kl}f_2(\hat{\bold R}^l,\hat{\bold \Theta}^l)\rangle\\
&+\Delta  t a_{kk}\left(\frac{1}{\varepsilon^2}(\bold B(\hat{\bold \Theta}^k)-\langle\hat{\bold R}^k\rangle)+\frac{\mu}{3}\partial_{xx}\bold B(\hat{\bold \Theta}^k)\right).
\end{aligned}
\end{equation}

Add (\ref{64}) and (\ref{65}),
\begin{equation}\label{66}
\begin{aligned}
\hat{\boldsymbol\Theta}^k+\langle \hat{\bold R}^k\rangle=&\hat{\boldsymbol \theta}^n+\Delta t\sum_{l=1}^{k-1}a_{kl}h(\hat{\boldsymbol\Theta}^l)\\
&+\langle \hat{\bold r}^n+\Delta t\sum_{l=1}^{k-1}\tilde{a}_{kl}f_1(\hat{\bold J}^l,\hat{\bold \Theta}^l)+\Delta t\sum_{l=1}^{k-1}a_{kl}f_2(\hat{\bold R}^l,\hat{\boldsymbol \Theta}^l)\rangle\\
&+\Delta ta_{kk}\left(\partial_{xx}\hat{\boldsymbol \Theta}^k+\frac{\mu}{3}\partial_{xx}\bold B(\hat{\boldsymbol \Theta}^k)\right).
\end{aligned}
\end{equation}

Thus, $\langle\hat{\bold  R}^k\rangle$ can be expressed in terms of $\hat{\boldsymbol\Theta}^k$ and other explicit terms. Inserting it back to (\ref{64}), then
\begin{equation}\label{67}
\begin{aligned}
\hat{\boldsymbol \Theta}^k=&\hat{\boldsymbol \theta}^n+\Delta t\sum_{l=1}^{k-1}a_{kl}h(\hat{\boldsymbol \Theta}^l)+\Delta ta_{kk}\partial_{xx}\hat{\boldsymbol \Theta}^k-\frac{\Delta ta_{kk}}{\varepsilon^2}\bold B(\hat{\boldsymbol \Theta}^k)\\
&+\frac{\Delta ta_{kk}}{\varepsilon^2}\left[\hat{\boldsymbol \theta}^n+\Delta t\sum_{l=1}^{k-1}a_{kl}h(\hat{\boldsymbol \Theta}^l)+\langle \hat{\bold r}^n+\Delta t\sum_{l=1}^{k-1}\tilde{a}_{kl}f_1(\hat{\bold J}^l,\hat{\bold \Theta}^l)+\Delta t\sum_{l=1}^{k-1}a_{kl}f_2(\hat{\bold R}^l,\hat{\boldsymbol \Theta}^l)\rangle\right.\\
&\left.+\Delta ta_{kk}\left(\partial_{xx}\hat{\boldsymbol \Theta}^k+\frac{\mu}{3}\partial_{xx}\bold B(\hat{\boldsymbol \Theta}^k)\right)-\hat{\boldsymbol \Theta}^k\right].
\end{aligned}
\end{equation}

Reorder and approximate $\bold B(\hat{\boldsymbol \Theta}^k)$ by $\bold B(\hat{\boldsymbol\Theta}^{k-1})+4\bold C(\hat{\boldsymbol\Theta}^{k-1})(\hat{\boldsymbol \Theta}^k-\hat{\boldsymbol\Theta}^{k-1})$ and $\partial_{xx}\bold B(\hat{\boldsymbol \Theta}^k)$ by $\partial_x(4\bold C(\hat{\boldsymbol \Theta}^{k-1})\partial_x\hat{\boldsymbol \Theta}^k)$, one gets

\begin{equation}\label{68}
\begin{aligned}
&\left[\left(1+\frac{\Delta t a_{kk}}{\varepsilon^2}\right)\bold I+\frac{\Delta t a_{kk}}{\varepsilon^2}4\bold C(\hat{\boldsymbol\Theta}^{k-1})\right]\hat{\bold \Theta}^k-\Delta t a_{kk}\left(1+\frac{\Delta t a_{kk}}{\varepsilon^2}\right)\partial_{xx}\hat{\boldsymbol\Theta}^k\\
&-\frac{(\Delta t a_{kk})^2}{\varepsilon^2}\frac{4\mu}{3}\partial_x(\bold C(\hat{\boldsymbol \Theta}^{k-1})\partial_x\hat{\boldsymbol \Theta}^k)\\
=&\left(1+\frac{\Delta t a_{kk}}{\varepsilon^2}\right)\hat{\boldsymbol \theta}^n+\left(1+\frac{\Delta t a_{kk}}{\varepsilon^2}\right)\Delta t\sum_{l=1}^{k-1}a_{kl}h(\hat{\boldsymbol\Theta}^l)\\
&+\frac{\Delta ta_{kk}}{\varepsilon^2}\left[\langle \hat{\bold r}^n+\Delta t\sum_{l=1}^{k-1}\tilde{a}_{kl}f_1(\hat{\bold J}^l,\hat{\boldsymbol \Theta}^l)+\Delta t\sum_{l=1}^{k-1}a_{kl}f_2(\hat{\bold R}^l,\hat{\boldsymbol \Theta}^l)\rangle+3\bold B(\hat{\boldsymbol\Theta}^{k-1})\right],
\end{aligned}
\end{equation}
where $\bold I$ is a $K\times K$ identity matrix, $\bold B$ and $\bold C$ are the same as (\ref{56}) and (\ref{58}).

Following the same space discretization as in section 2.3,
\begin{equation}\label{69}
\begin{aligned}
&\left[\left(1+\frac{\Delta t a_{kk}}{\varepsilon^2}\right)\bold I+\frac{\Delta t a_{kk}}{\varepsilon^2}4\bold C(\hat{\boldsymbol\Theta})_i^{k}\right]\hat{\bold \Theta}_i^k-\Delta t a_{kk}\left(1+\frac{\Delta t a_{kk}}{\varepsilon^2}\right)\frac{\hat{\boldsymbol\Theta}_{i+1}^k-2\hat{\boldsymbol\Theta}_{i}^k+\hat{\boldsymbol\Theta}_{i-1}^k}{\Delta x^2}\\
&-\frac{(\Delta t a_{kk})^2}{\varepsilon^2}\frac{4\mu}{3}\frac{\bold C(\hat{\boldsymbol \Theta})_{i+1/2}^{k-1}(\hat{\boldsymbol \Theta}_{i+1}^k-\hat{\boldsymbol \Theta}_{i}^k)-\bold C(\hat{\boldsymbol \Theta})_{i-1/2}^{k-1}(\hat{\boldsymbol \Theta}_{i}^k-\hat{\boldsymbol \Theta}_{i-1}^k)}{\Delta x^2}\\
=&\left(1+\frac{\Delta t a_{kk}}{\varepsilon^2}\right)\hat{\boldsymbol \theta}_i^n+\left(1+\frac{\Delta t a_{kk}}{\varepsilon^2}\right)\Delta t\sum_{l=1}^{k-1}a_{kl}h(\hat{\boldsymbol \Theta}_i^l)\\
&+\frac{\Delta ta_{kk}}{\varepsilon^2}\left[\langle \hat{\bold r}_i^n+\Delta t\sum_{l=1}^{k-1}\tilde{a}_{kl}f_1(\hat{\bold J}_i^l,\hat{\boldsymbol \Theta}_i^l)+\Delta t\sum_{l=1}^{k-1}a_{kl}f_2(\hat{\bold R}_i^l,\hat{\boldsymbol \Theta}_i^l)\rangle+3\bold B(\hat{\boldsymbol\Theta}_i^{k-1})\right],
\end{aligned}
\end{equation}
where 
$$
\begin{aligned}
&h(\hat{\boldsymbol \Theta}_i^l)=\frac{\hat{\boldsymbol \Theta}_{i+1}^l-2\hat{\boldsymbol \Theta}_i^l+\hat{\boldsymbol \Theta}_{i-1}^l}{\Delta x^2}-\frac{1}{\varepsilon^2}(\bold B (\hat{\boldsymbol \Theta}_i^l)-\langle \hat{\bold R}_i^l\rangle),\\
&f_1(\hat{\bold J}_i^l, \hat{\boldsymbol \Theta}_i^l)=-\frac{v}{2\Delta x}(\hat{\bold J}_{i+1}^l-\hat{\bold J}_{i-1}^l)+\frac{v\phi^{1/2}}{2\Delta x}(\hat{\bold R}_{i+1}^l-2\hat{\bold R}_i^l+\hat{\bold R}_{i-1}^l)\\
& \ \ \ \ \ \ \ \ \ \ \ \ \ \ \ \ -\frac{v\phi^{1/2}}{4}(\boldsymbol\gamma_i^l-\boldsymbol\gamma_{i-1}^l+\boldsymbol\beta_{i+1}^l-\boldsymbol \beta_i^l)-\frac{\mu}{3}\frac{\bold B(\hat{\boldsymbol \Theta}_{i+1}^l)-2\bold B(\hat{\boldsymbol \Theta}_i^l)+\bold B(\hat{\boldsymbol \Theta}_{i-1}^l)}{\Delta x^2},\\
&f_2(\hat{\bold R}_i^l,\hat{\boldsymbol \Theta}_i^l)=\frac{1}{\varepsilon^2}(\bold B(\hat{\boldsymbol \Theta}_i^l)-\hat{\bold R}_i^l)+\frac{\mu}{3}\frac{\bold B(\hat{\boldsymbol \Theta}_{i+1}^l)-2\bold B(\hat{\boldsymbol \Theta}_{i}^l)+\bold B(\hat{\boldsymbol \Theta}_{i-1}^l)}{\Delta x^2},
\end{aligned}
$$
with $\boldsymbol\gamma^l_i$ and $\boldsymbol\beta^l_i$ the same as (\ref{27}).

We omit the details of velocity discretization and other space discretization, which is the same as in section 2.3. On the left hand side of (\ref{69}), we specify the approximation of $\bold C(\hat{\boldsymbol\Theta})_{i+1/2}^{k-1}$ by 
$$\bold C(\hat{\boldsymbol\Theta})_{i+1/2}^{k-1}\approx\frac{1}{2}\left(\bold C(\hat{\boldsymbol\Theta})_{i+1}^{k-1}+\bold C(\hat{\boldsymbol\Theta})_{i}^{k-1}\right).$$

Notice that to solve out $\hat{\boldsymbol \Theta}_i^k$, one needs to invert a large $K(N_x+1)\times K(N_x+1)$ matrix $\bold M$, where $N_x+1$ is the number of discretized points in space. $\bold M$ is defined as following:
\begin{equation}\label{71}
\bold M_{Ki+1:K(i+1),Kj+1:K(j+1)}=\left\{
\begin{aligned}
&\bold m^1_{i}&&\mbox{for} \ j=i-1\\
&\bold m^2_i&&\mbox{for} \ j=i\\
&\bold m^3_{i}&&\mbox{for} \ j=i+1\\
&0&&\mbox{otherwise}
\end{aligned}, \ \ \ i=0,1,\cdots, N_x
\right.
\end{equation}
with $K\times K$ matrix $\bold m^1, \bold m^2$ and $\bold m^3$:
\begin{equation}\label{72}
\begin{aligned}
\bold m^1_i=&-\frac{\Delta t a_{kk}}{\Delta x^2}\left(1+\frac{\Delta t a_{kk}}{\varepsilon^2}\right)\bold I-\frac{(\Delta ta_{kk})^2}{2\varepsilon^2\Delta x^2}\frac{4\mu}{3}\left[\bold C(\hat{\boldsymbol \Theta})_{i-1}^{k-1}+\bold C (\hat{\boldsymbol \Theta})_i^{k-1}\right];\\
\bold m^2_i=&\left(1+\frac{\Delta t a_{kk}}{\varepsilon^2}\right)\bold I+\frac{\Delta t a_{kk}}{\varepsilon^2}4\bold C(\hat{\boldsymbol\Theta})_i^{k-1}+\frac{2\Delta t a_{kk}}{\Delta x^2}\left(1+\frac{\Delta t a_{kk}}{\varepsilon^2}\right)\bold I\\
&+\frac{(\Delta ta_{kk})^2}{2\varepsilon^2\Delta x^2}\frac{4\mu}{3}\left[\bold C(\hat{\boldsymbol \Theta})_{i+1}^{k-1}+2\bold C (\hat{\boldsymbol \Theta})_i^{k-1}+\bold C (\hat{\boldsymbol \Theta})_{i-1}^{k-1}\right];\\
\bold m^3_i=&-\frac{\Delta t a_{kk}}{\Delta x^2}\left(1+\frac{\Delta t a_{kk}}{\varepsilon^2}\right)\bold I-\frac{(\Delta ta_{kk})^2}{2\varepsilon^2\Delta x^2}\frac{4\mu}{3}\left[\bold C(\hat{\boldsymbol \Theta})_{i}^{k-1}+\bold C (\hat{\boldsymbol \Theta})_{i+1}^{k-1}\right].\end{aligned}
\end{equation}

From \cite{JinLu}, one can prove the matrix $\left(1+\frac{\Delta t a_{kk}}{\varepsilon^2}\right)\bold I+\frac{\Delta t a_{kk}}{\varepsilon^2}4\bold C(\hat{\boldsymbol\Theta})_i^{k-1}$ is symmetric, positive and definite under some gentle conditions. Then obviously, $\bold M$ is strictly block diagonal dominant, thus invertible. Since $\bold M$ is also sparse, fast algorithm can be applied to solve $\hat{\boldsymbol \Theta}^k$, then $\hat{\bold R}^k$ and $\hat{\bold J}^k$ can be obtained from (\ref{63}a) and (\ref{63}b) subsequently.

\subsection{ The sAP Property}
Similar as in section 2.4, setting $\varepsilon\to0$ in (\ref{63}), one gets
\begin{subequations}\label{73}
\begin{align}
&\hat{\bold R}_i^k=\bold B(\hat{\boldsymbol\Theta}_i^k),\\
&\hat{\bold J}_i^k=-v\frac{\hat{\bold R}_{i+1}^k-\hat{\bold R}_{i-1}^k}{2\Delta x},\\
&\langle \hat{\bold R}_i^k\rangle=\bold B(\hat{\boldsymbol \Theta}_i^k).
\end{align}
\end{subequations}

Inserting (\ref{73}) back into (\ref{62}a) $+$ (\ref{62}c) and setting $\varepsilon=0$,
\begin{equation}\label{74}
\hat{\bold r}_i^{n+1}+\hat{\boldsymbol \theta}_i^{n+1}= \hat{\bold r}_i^n+\hat{\boldsymbol\theta}_i^n+\Delta t\sum_{k=1}^s\tilde{b}_k \hat f_1(\hat{\bold R}_i^k,\hat{\boldsymbol\Theta}_i^k)+\Delta t\sum_{k=1}^sb_k\hat f_2(\hat{\bold R}_i^k,\hat{\boldsymbol\Theta}_i^k)+\Delta t\sum_{k=1}^sb_k\hat h(\hat{\boldsymbol\Theta}_i^k),
\end{equation}
where
\begin{subequations}\label{75}
\begin{align}
&\hat f_1(\hat{\bold R}^k,\hat{\boldsymbol\Theta}^k)=v^2\frac{\hat{\bold R}_{i+2}^k-2\hat{\bold R}_{i}^k+\hat{\bold R}_{i-2}^k}{4\Delta x^2}-\frac{\bold B(\hat{\boldsymbol\Theta}_{i+1}^k)-2\bold B(\hat{\boldsymbol\Theta}_{i}^k)+\bold B(\hat{\boldsymbol\Theta}_{i-1}^k)}{3\Delta x^2},\\
&\hat f_2(\hat{\boldsymbol\Theta}^k)=\frac{\bold B(\hat{\boldsymbol\Theta}_{i+1}^k)-2\bold B(\hat{\boldsymbol\Theta}_{i}^k)+\bold B(\hat{\boldsymbol\Theta}_{i-1}^k)}{3\Delta x^2},\\
&\hat h(\hat{\boldsymbol\Theta}^k)=\frac{\hat{\boldsymbol\Theta}_{i+1}^k-2\hat{\boldsymbol\Theta}_{i}^k+\hat{\boldsymbol\Theta}_{i-1}^k}{\Delta x^2}.
\end{align}
\end{subequations}

Take $\langle\cdot \rangle$ on (\ref{74}) and use (\ref{73}c),
\begin{equation}\label{76}
\begin{aligned}
\hat{\boldsymbol \theta}_i^{n+1}+\bold B(\hat{\boldsymbol\theta}_i^{n+1})=&\hat{\boldsymbol \theta}_i^n+\bold B(\hat{\boldsymbol \theta}_i^n)+\Delta t\sum_{k=1}^sb_k\left[\frac{\hat{\boldsymbol\Theta}_{i+1}^k-2\hat{\boldsymbol\Theta}_{i}^k+\hat{\boldsymbol\Theta}_{i-1}^k}{\Delta x^2}\right.\\
&\left.+\frac{\bold B(\hat{\boldsymbol\Theta}_{i+1}^k)-2\bold B(\hat{\boldsymbol\Theta}_{i}^k)+\bold B(\hat{\boldsymbol\Theta}_{i-1}^k)}{3\Delta x^2}\right]+O(\Delta x^2),
\end{aligned}
\end{equation}
where
\begin{equation}\label{77}
\begin{aligned}
\hat{\boldsymbol \Theta}_i^{k}+\bold B(\hat{\boldsymbol\Theta}_i^{k})=&\hat{\boldsymbol \theta}_i^n+\bold B(\hat{\boldsymbol \theta}_i^n)+\Delta t\sum_{l=1}^ka_{kl}\left[\frac{\hat{\boldsymbol\Theta}_{i+1}^l-2\hat{\boldsymbol\Theta}_{i}^l+\hat{\boldsymbol\Theta}_{i-1}^l}{\Delta x^2}\right.\\
&\left.+\frac{\bold B(\hat{\boldsymbol\Theta}_{i+1}^l)-2\bold B(\hat{\boldsymbol\Theta}_{i}^l)+\bold B(\hat{\boldsymbol\Theta}_{i-1}^l)}{3\Delta x^2}\right],
\end{aligned}
\end{equation}
which is an {\it implicit} R-K method for the limiting diffusion equation (\ref{57}). Thus the sAP property is illustrated.

\section{Numerical Tests}
In this section we present several numerical tests to illustrate the performance of the new sAP IMEX schemes. We use notations SSP2 to denote the new scheme based on the second order IMEX-RK tableu (\ref{eq:SSP2}) and JPT to denote the second order scheme in \cite{jin2000uniformly}. First we test the new approach in the deterministic case and then we consider the presence of uncertainties.

\subsection{Test 1: The Linear Transport Equation}
We consider problem (\ref{1}) with 
$$x\in[0,1], \ F_L(v)=1, \ F_R(v)=0, \ \sigma_s=1,\  \sigma_A=0, \ Q=0,\ \varepsilon=10^{-6}$$
We reported the results at different time $t=0.01,t=0.05$ and $t=0.15$ in Figure 1. In all figures we use notations $NT$ to denote the number of time iterations. With fixed $\Delta x=0.01$, for JPT ($\circ$) $\Delta t=0.5(\Delta x)^2$, $NT=200,1000,3000$; for SSP2 ($\diamond$) $\Delta t=\lambda \Delta x$ with $\lambda=0.04$, $NT=25,125,375$. The exact diffusive solution is computed by equation (2) with the same mesh size in dashed line. One can see that both methods are AP and work well in the diffusive regime at any times. The time cost of the two methods are reported in Table 1, which shows that SSP2 is about 7 times more efficient than JPT.

\begin{table}[h]
\begin{center}
\begin{tabular}{ c|c|c|c } 
 \hline
  & $t=0.01$& $t=0.05$&$t=0.15$ \\ 
  \hline
 JPT & 0.451567s&2.014770s&6.069805s \\ 
  \hline
SSP2 &0.069801s &0.290240s& 0.815041s\\ 
 \hline
\end{tabular}
\end{center}
\caption{Test 1. Time cost comparison in seconds between JPT and SSP2 schemes.}
\end{table}

\begin{figure}[tb]
\begin{center}
\includegraphics[scale=0.45]{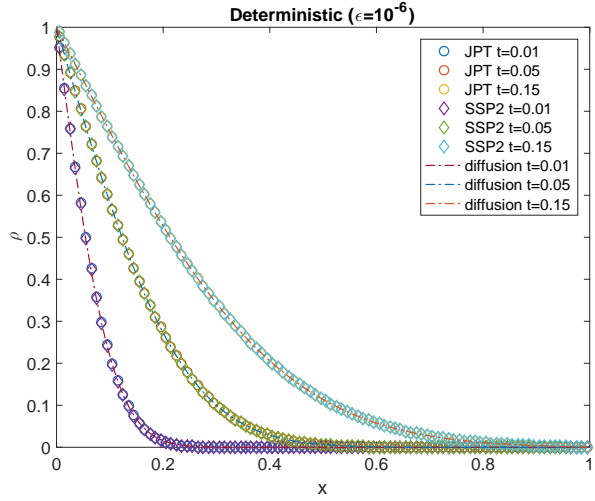}
\end{center}
\caption{Test 1. The solution of the mass density $\rho$ to the linear transport equation. JPT ($\circ$), $\Delta x=0.01$, $\Delta t=0.5(\Delta x)^2$, $NT=200, 1000, 3000$; SSP2 ($\diamond$), $\Delta x=0.01$, $\Delta t=\lambda\Delta x$ with $\lambda=0.04$, $NT=25, 125, 375$.}
\end{figure}

\subsection{Test 2: The Linear Transport Equation with Random Inputs}
We consider problem (\ref{1}) with random inputs and 
$$x\in[0,1], \ F_L(v)=1, \ F_R(v)=0, \ \sigma_s=1+0.5z,\  \sigma_A=0, \ Q=0,\ \varepsilon=10^{-6},$$
where $z$ follows uniform distribution on $[-1,1]$ (Denoted by $z\sim \mathcal U[-1,1]$ in the following tests).

The mean and standard deviation are two quantities of our interest. So we examine these two quantities of the numerical solutions at different time $t=0.01,t=0.05$ and $t=0.15$ in Figure 2. Given the gPC coefficients $\hat{\rho}_k$ of $\rho$, it's convenient to calculate the mean value and standard deviation by
$$\text{E}[\rho]\approx \hat \rho_1, \ \text{Sd}[\rho]\approx \sqrt{\sum_{k=2}^K\hat{\rho}_k^2}.$$
With fixed $\Delta x=0.025$, for JPT ($\circ$) $\Delta t=0.0002$ and $NT=50,250,750$; for SSP2 ($\diamond$) $\Delta t=\lambda \Delta x$ with $\lambda=0.035$, $NT=11,57,171$. The exact diffusive solution is represented by the dashed line. As expected, both methods are sAP and match well in the diffusive regime at any times. The distribution of $\rho(x,z)$ at different times are presented in Figure 3. We report the time cost comparison of gPC-JPT and gPC-SSP2 in Table 2, which shows that gPC-SSP2 is more than 3 times faster than gPC-JPT. Since randomness is introduced and computation becomes more time-consuming, the efficiency of gPC-SSP2 is remarkable.

\begin{table}[h!]
\begin{center}
\begin{tabular}{ c|c|c|c } 
 \hline
  & $t=0.01$& $t=0.05$&$t=0.15$ \\ 
  \hline
 gPC-JPT & 2.662382s&13.159598s&39.255826s \\ 
  \hline
gPC-SSP2 &0.759259s &3.865332s&11.424931s\\ 
 \hline
\end{tabular}
\end{center}
\caption{Test 2. Time cost comparison in seconds between JPT and SSP2 schemes.}
\end{table}

\begin{figure}[tb]
\includegraphics[scale=0.4]{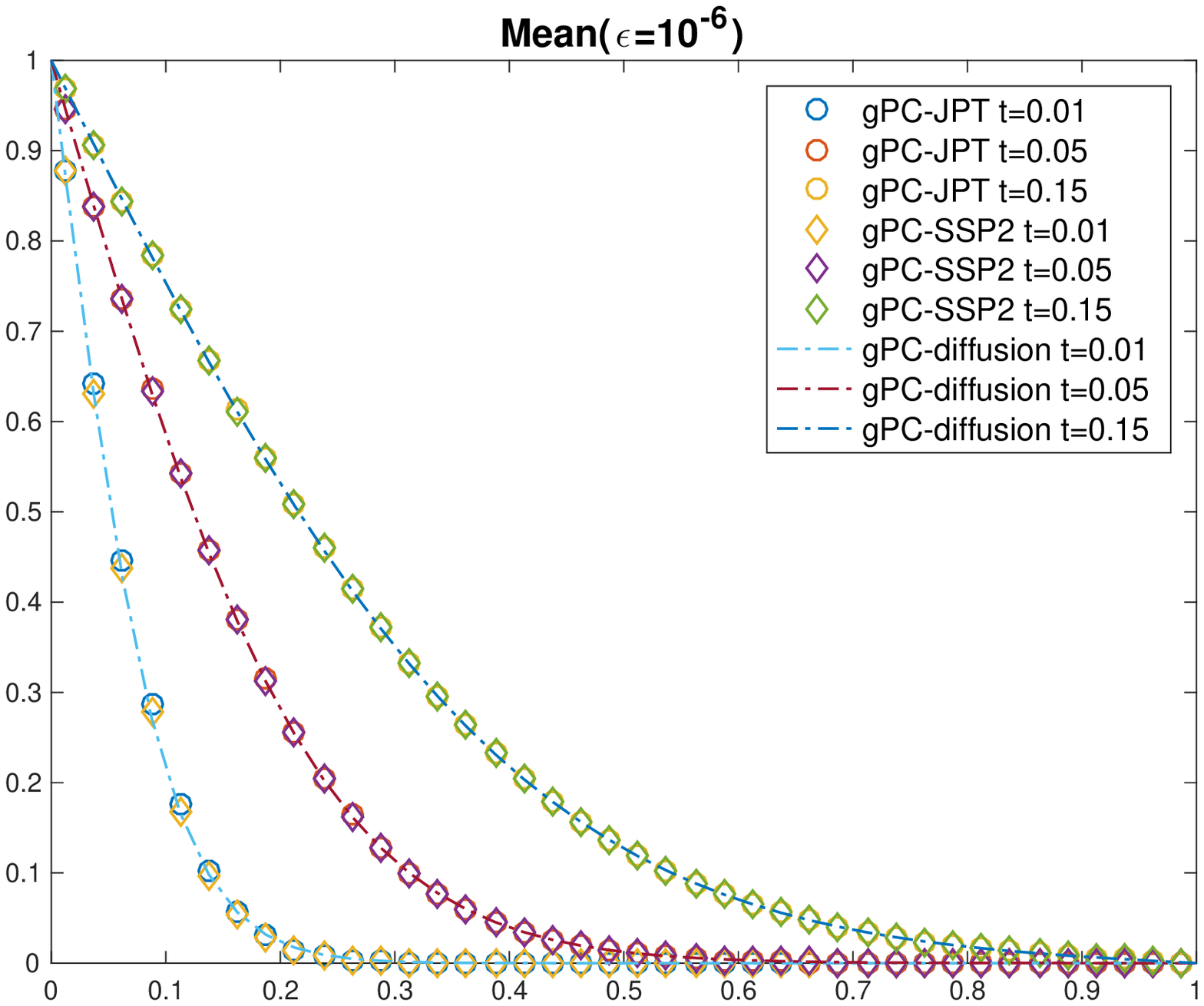}
\includegraphics[scale=0.4]{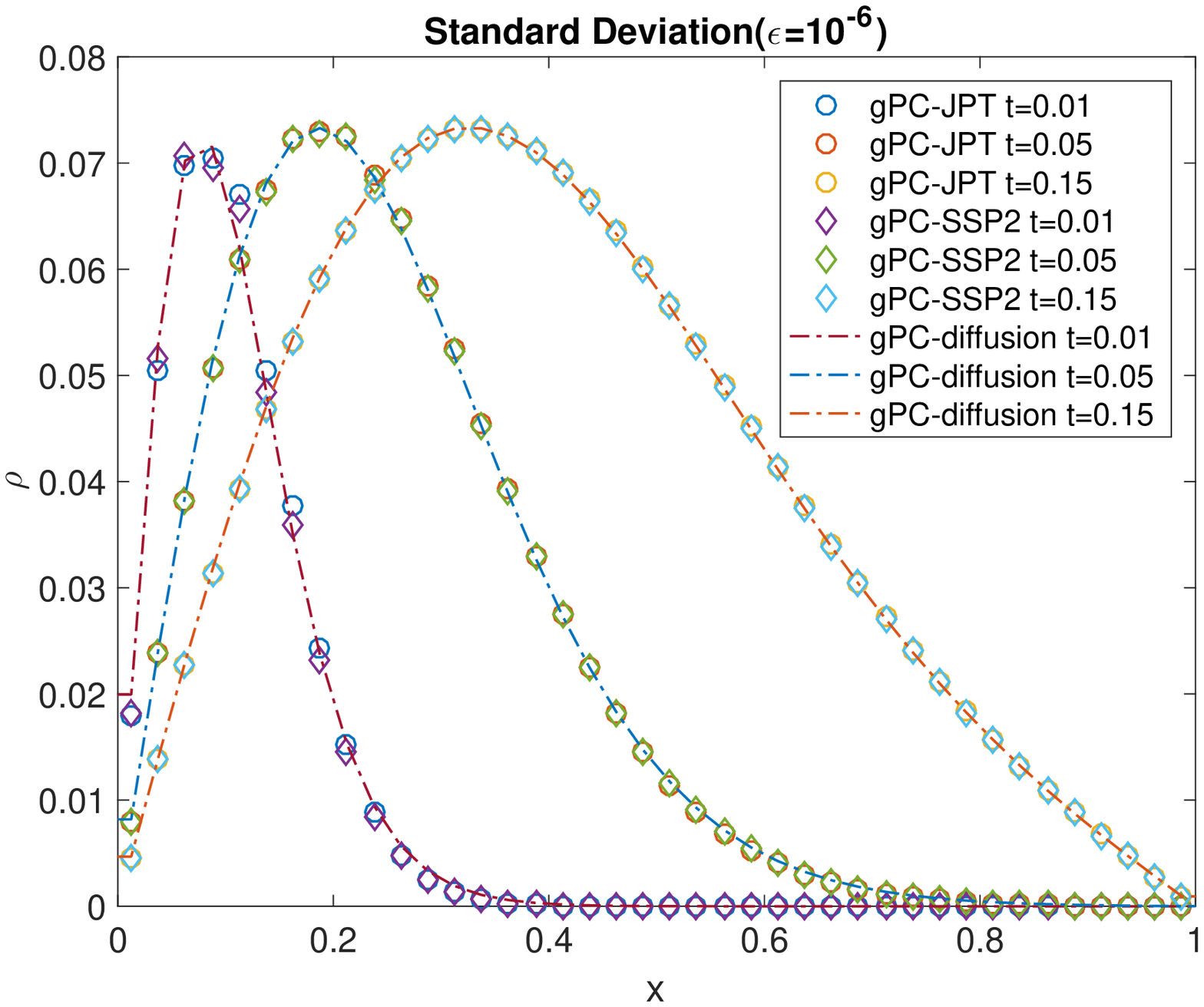}
\caption{Test 2. The mean (left) and standard deviation (right) of the solution of the mass density $\rho$ to the linear transport equation. JPT ($\circ$), $\Delta x=0.025$, $\Delta t=0.0002$, $NT=50, 250, 750$; SSP2 ($\diamond$), $\Delta x=0.025$, $\Delta t=\lambda\Delta x$ with $\lambda=0.035$, $NT=11, 57, 171$.}
\end{figure}

\begin{figure}[tb]
\includegraphics[scale=0.3]{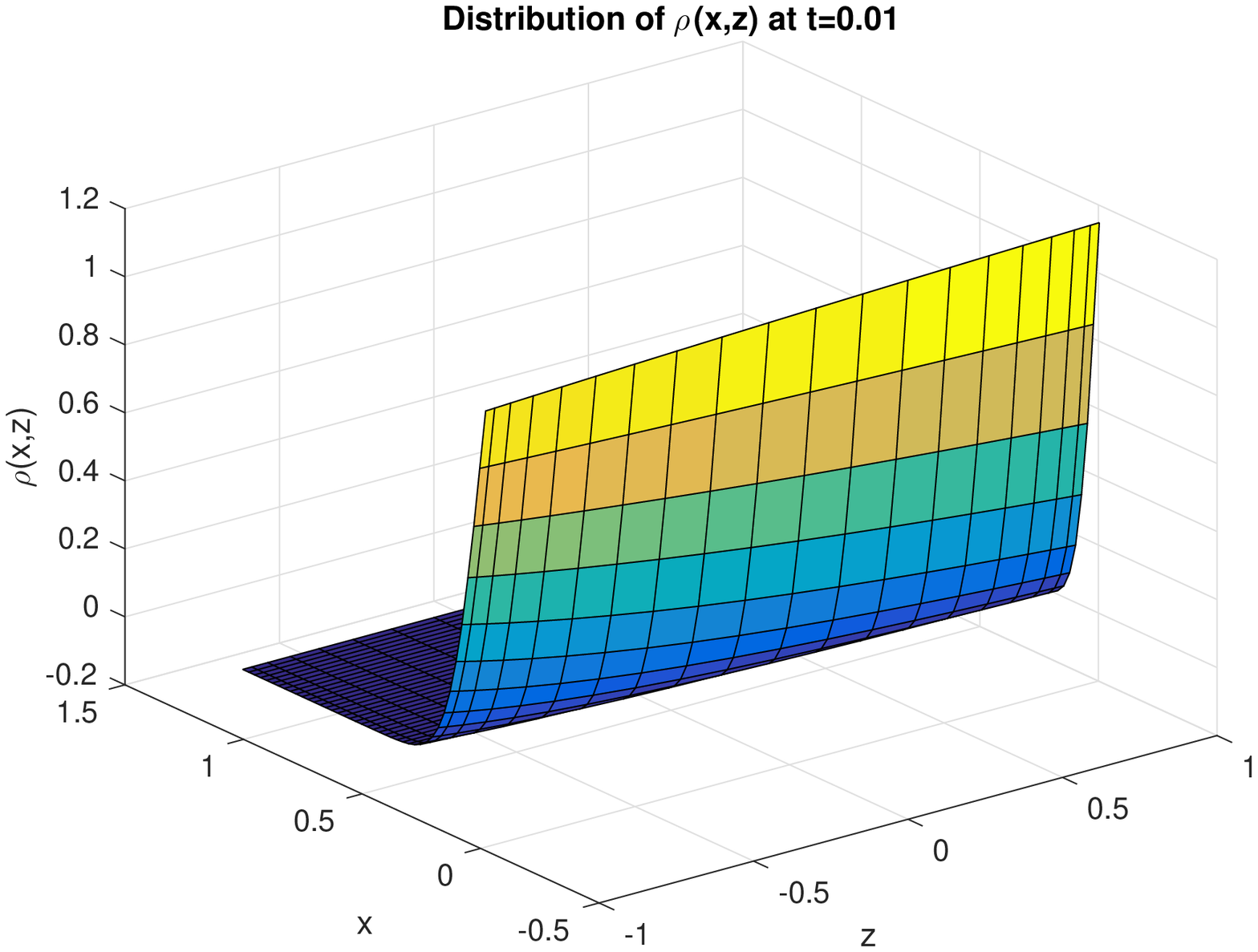}\includegraphics[scale=0.3]{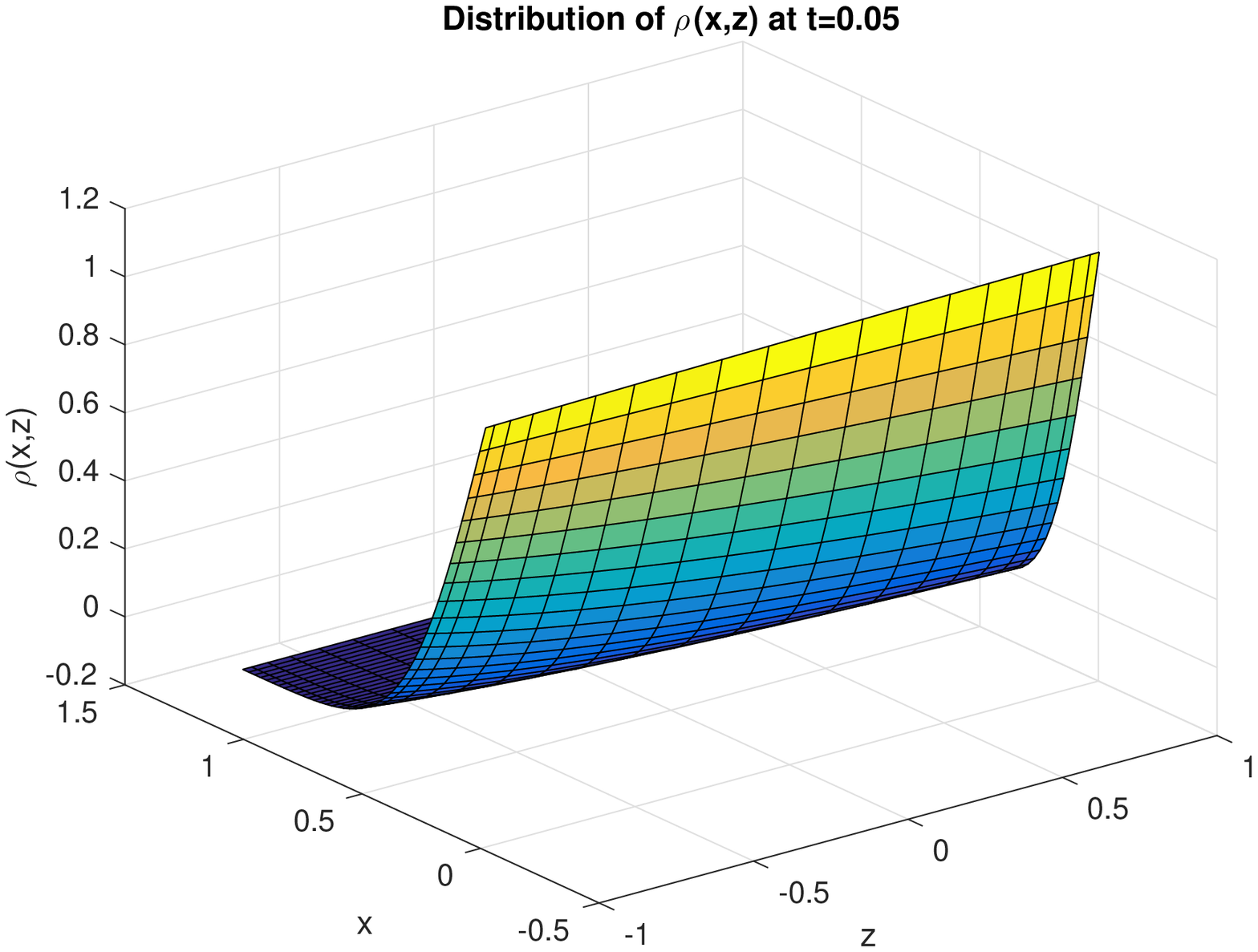}\includegraphics[scale=0.3]{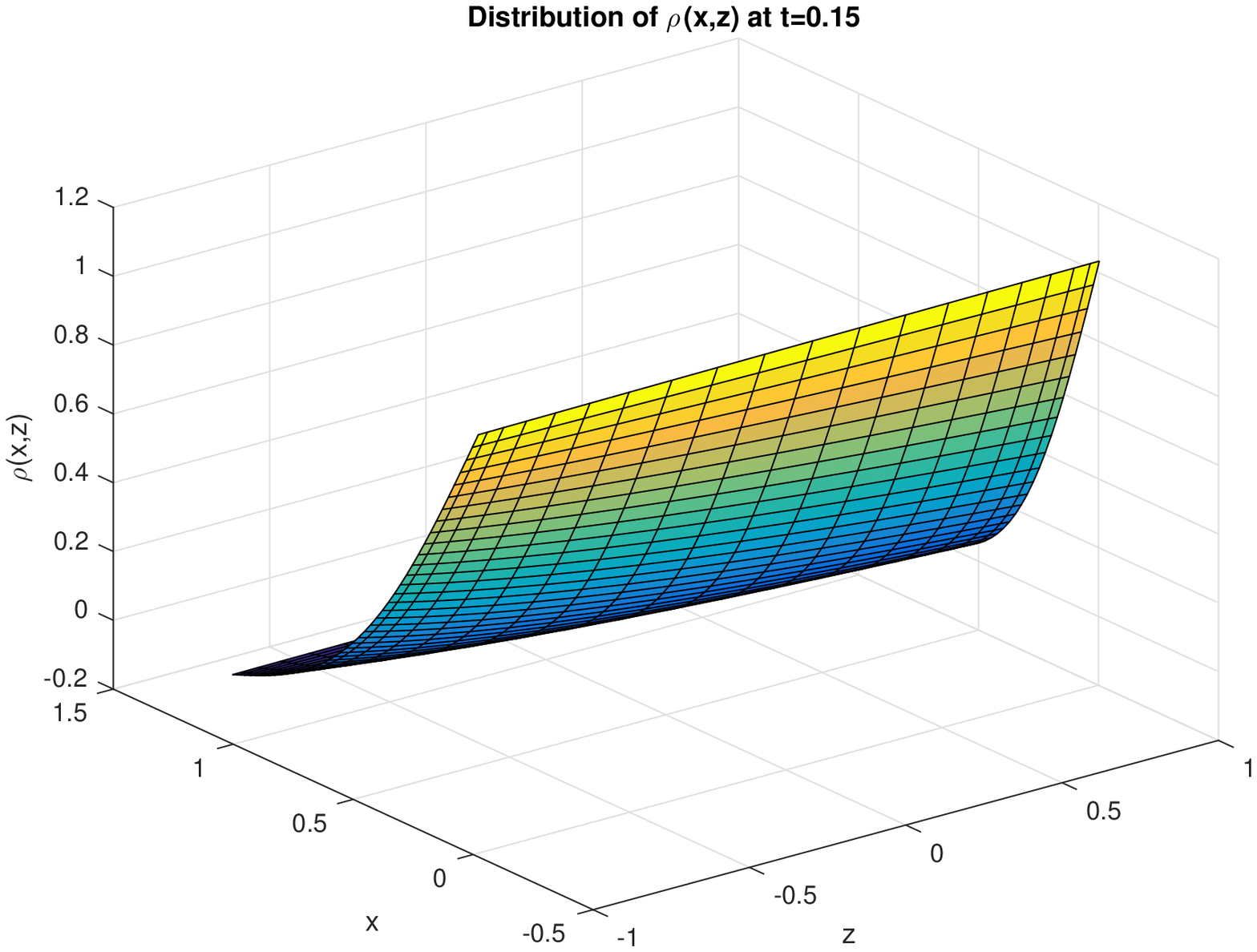}
\caption{Test 2. Distribution of the mass density $\rho(x,z)$ to the linear transport equation by gPC-SSP2 with $NT=11, 57, 171$.}
\end{figure}



\subsection{Test 3: The Radiative Heat Transfer Equation}
We then compare the two methods for deterministic radiative heat transfer equations (\ref{43}) near diffusive regime ($\varepsilon=10^{-6}$) with the following initial and boundary conditions:
$$
\begin{aligned}
&I_I(x,\mu,z,0)=0,\ \ \theta_I(x,z,0)=0, \ \ x\in[0,1]\\
&\theta_B(0,z,t)=1, \ \theta_B(1,z,t)=0; \\
&I_B(0,\mu,z,t)=1, \ \mu>0, \ \ I_B(1,\mu,z,t)=0, \ \mu<0,
\end{aligned}
$$
and constant coefficient
$$\sigma(x)=1.$$

Figure 4 shows the solutions of temperature $\theta$ at different time $t=0.01,t=0.05$ and $t=0.15$. With fixed $\Delta x=0.025$, for JPT ($\circ$) $\Delta t=0.002$; for SSP2 ($\diamond$) $\Delta t=\lambda \Delta x$ with $\lambda=0.035$. The exact diffusive solution is computed by equation (53) with the same mesh size in dashed line. One can see that both methods are AP and work well in the diffusive regime at any times. The time cost of both methods are reported in Table 3, which shows that SSP2 is about twice times faster than JPT in this nonlinear case.

\begin{table}[h!]
\begin{center}
\begin{tabular}{ c|c|c|c } 
 \hline
  & $t=0.01$& $t=0.05$&$t=0.15$ \\ 
  \hline
 JPT & 0.112507s&0.268415s&0.806410s \\ 
  \hline
SSP2 &0.059658s&0.208812s&0.583166s\\ 
 \hline
\end{tabular}
\end{center}
\caption{Test 3. Time cost comparison in seconds between JPT and SSP2 schemes.}
\end{table}

\begin{figure}[ht]
\begin{center}
\includegraphics[scale=0.45]{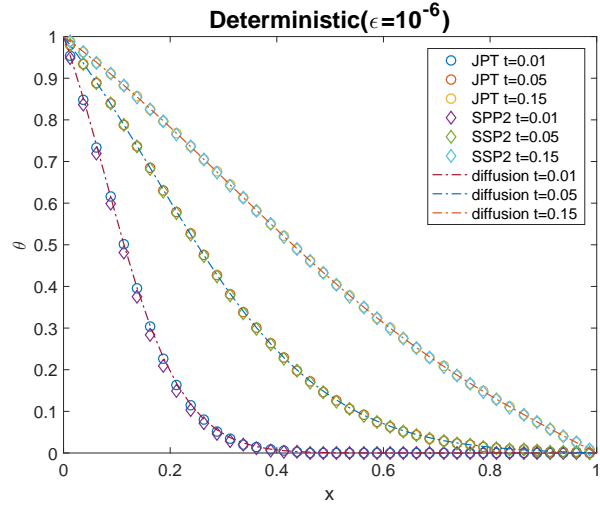}
\end{center}
\caption{Test 3. The solution of the temperature $\theta$ to the radiative heat transfer equation. JPT ($\circ$), $\Delta x=0.025$, $\Delta t=0.0002$, $NT=50, 250, 750$; SSP2 ($\diamond$), $\Delta x=0.025$, $\Delta t=\lambda\Delta x$ with $\lambda=0.035$, $NT=11, 57, 171$.}
\end{figure}

\subsection{Test 4: The Radiative Heat Transfer Equation with Random Inputs}
We then compare the two methods in the equations (\ref{43}) with randomness in the cross-section with the following initial and boundary conditions:
$$
\begin{aligned}
&I_I(x,\mu,z,0)=0,\ \ \theta_I(x,z,0)=0, \ \ x\in[0,1]\\
&\theta_B(0,z,t)=1, \ \theta_B(1,z,t)=0; \\
&I_B(0,\mu,z,t)=1+0.5z, \ \mu>0, \ \ I_B(1,\mu,z,t)=0, \ \mu<0,
\end{aligned}
$$
and random coefficient
$$\sigma(z)=1+0.5z, \ z\sim \mathcal U[-1,1].$$

In Figure 5, we report the mean and standard deviation of the temperature $\theta$ at different time $t=0.01,t=0.05$ and $t=0.15$. Similar to Test 2, the mean value and standard deviation are calculated by
$$\text{E}[\theta]\approx \hat \theta_1, \ \text{Sd}[\theta]\approx \sqrt{\sum_{k=2}^K\hat{\theta}_k^2}.$$
With fixed $\Delta x=0.025$, for JPT ($\circ$) $\Delta t=0.0002$ and $NT=50,250,750$; for SSP2 ($\diamond$) $\Delta t=\lambda \Delta x$ with $\lambda=0.035$, $NT=11,57,171$. The reference diffusive solution is represented by the dashed line. We refer to \cite{JinLu} for the same results using micro-macro decomposition based gPC method. Good agreements of gPC-JPT, gPC-SSP2 methods and reference solution are observed in the diffusive regime at any time. The distribution of $\theta(x,z)$ at different times are presented in Figure 6. We report the time cost comparison of gPC-JPT and gPC-SSP2 in Table 4, which shows that gPC-SSP2 is nearly 4 times faster than gPC-JPT. The nonlinearity slows down gPC-JPT a lot and thus the advantage in efficiency of gPC-SSP2 is significant in this nonlinear problem with random inputs.

\begin{table}[h!]
\begin{center}
\begin{tabular}{ c|c|c|c } 
 \hline
  & $t=0.01$& $t=0.05$&$t=0.15$ \\ 
  \hline
 gPC-JPT &172.431508s&804.915170s&2329.405654s \\ 
  \hline
gPC-SSP2 &36.945166s&192.040137s&626.736282s\\ 
 \hline
\end{tabular}
\end{center}
\caption{Test 4. Time cost comparison in seconds between JPT and SSP2 schemes.}
\end{table}

\begin{figure}[tb]
\includegraphics[scale=0.4]{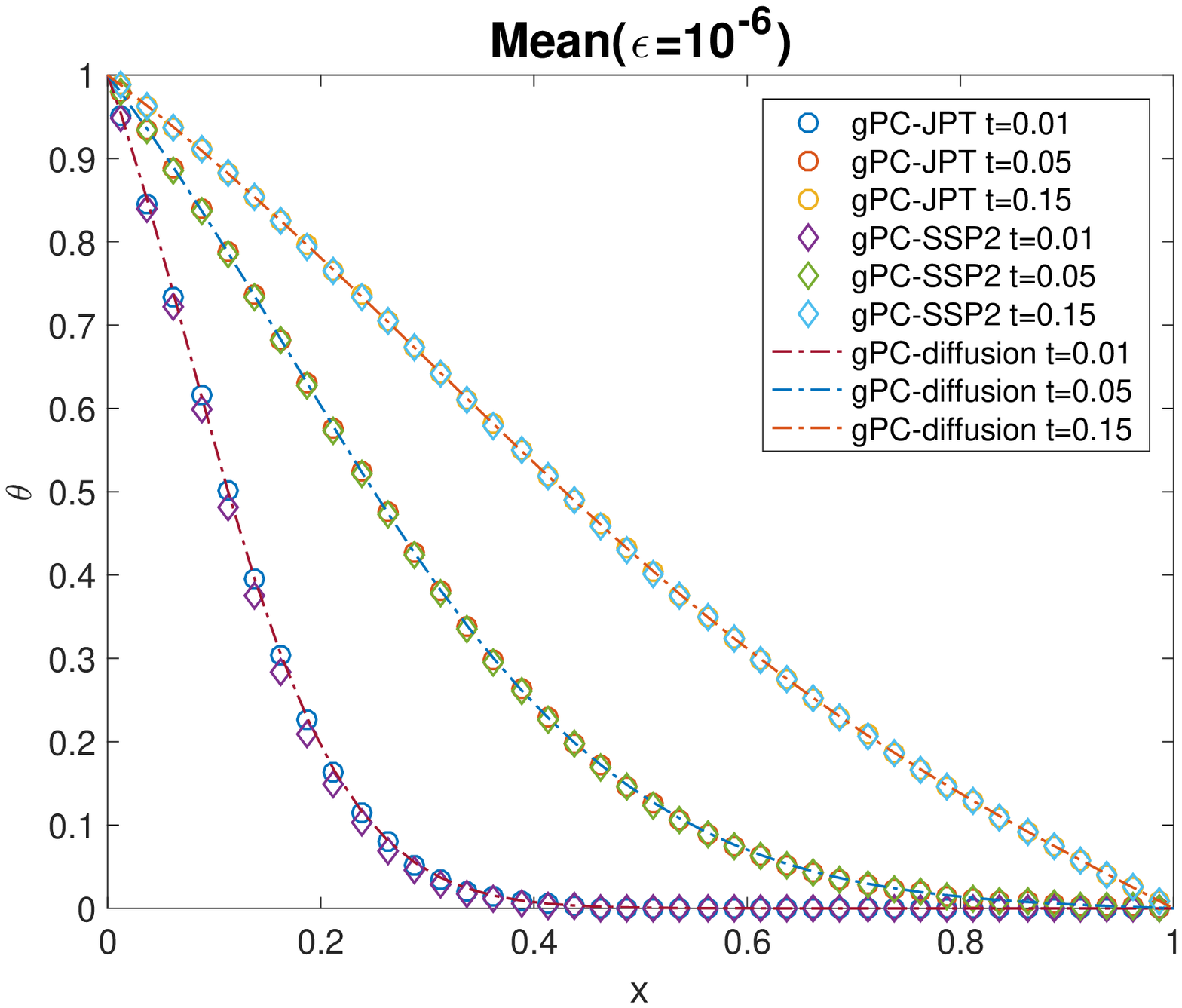}
\includegraphics[scale=0.4]{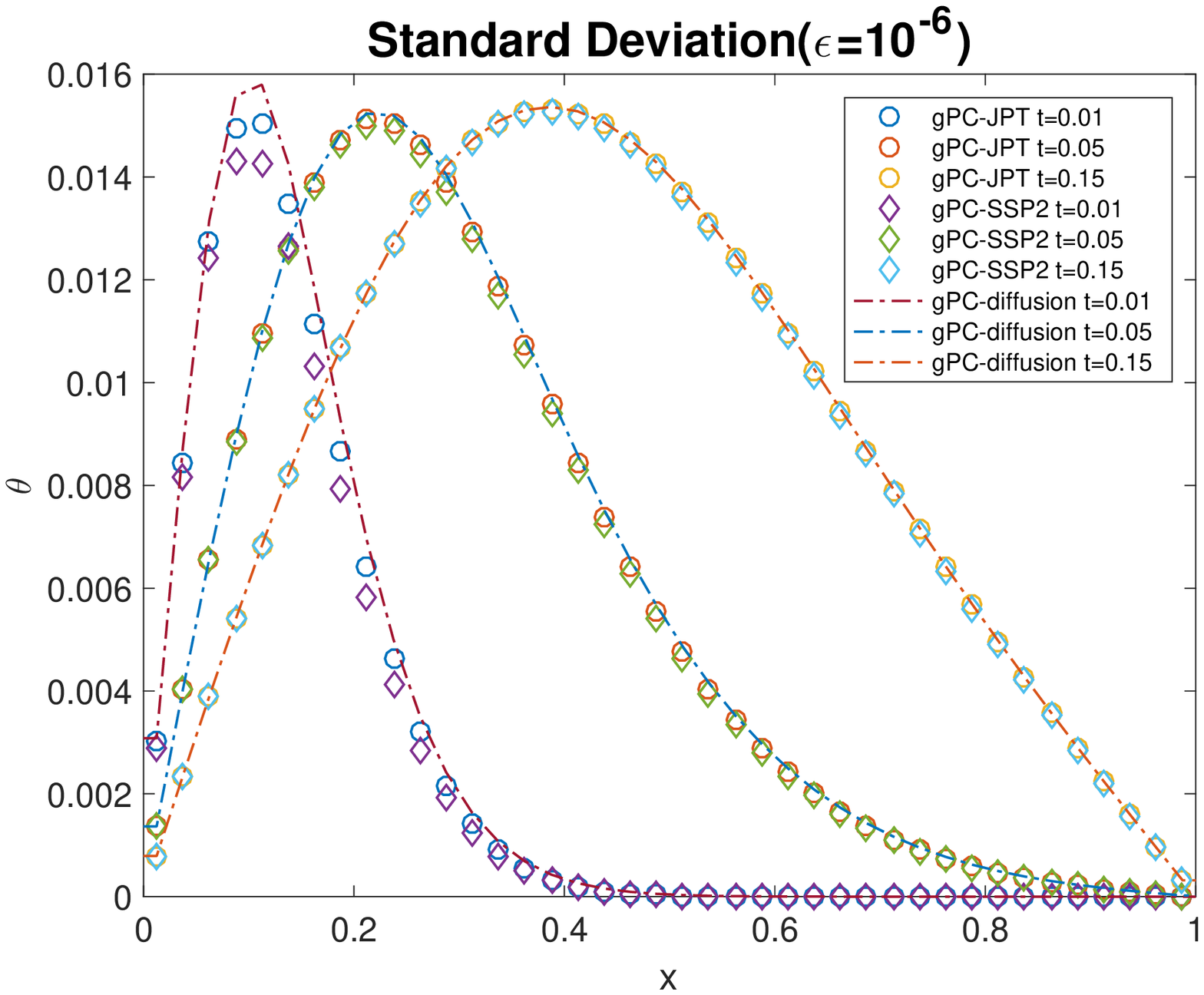}
\caption{Test 4. The mean (left) and standard deviation (right) of the solution of the temperature $\theta$ to the radiative heat transfer equation. JPT ($\circ$), $\Delta x=0.025$, $\Delta t=0.0002$, $NT=50, 250, 750$; SSP2 ($\diamond$), $\Delta x=0.025$, $\Delta t=\lambda\Delta x$ with $\lambda=0.035$, $NT=11, 57, 171$.}
\end{figure}
\begin{figure}[tb]
\includegraphics[scale=0.3]{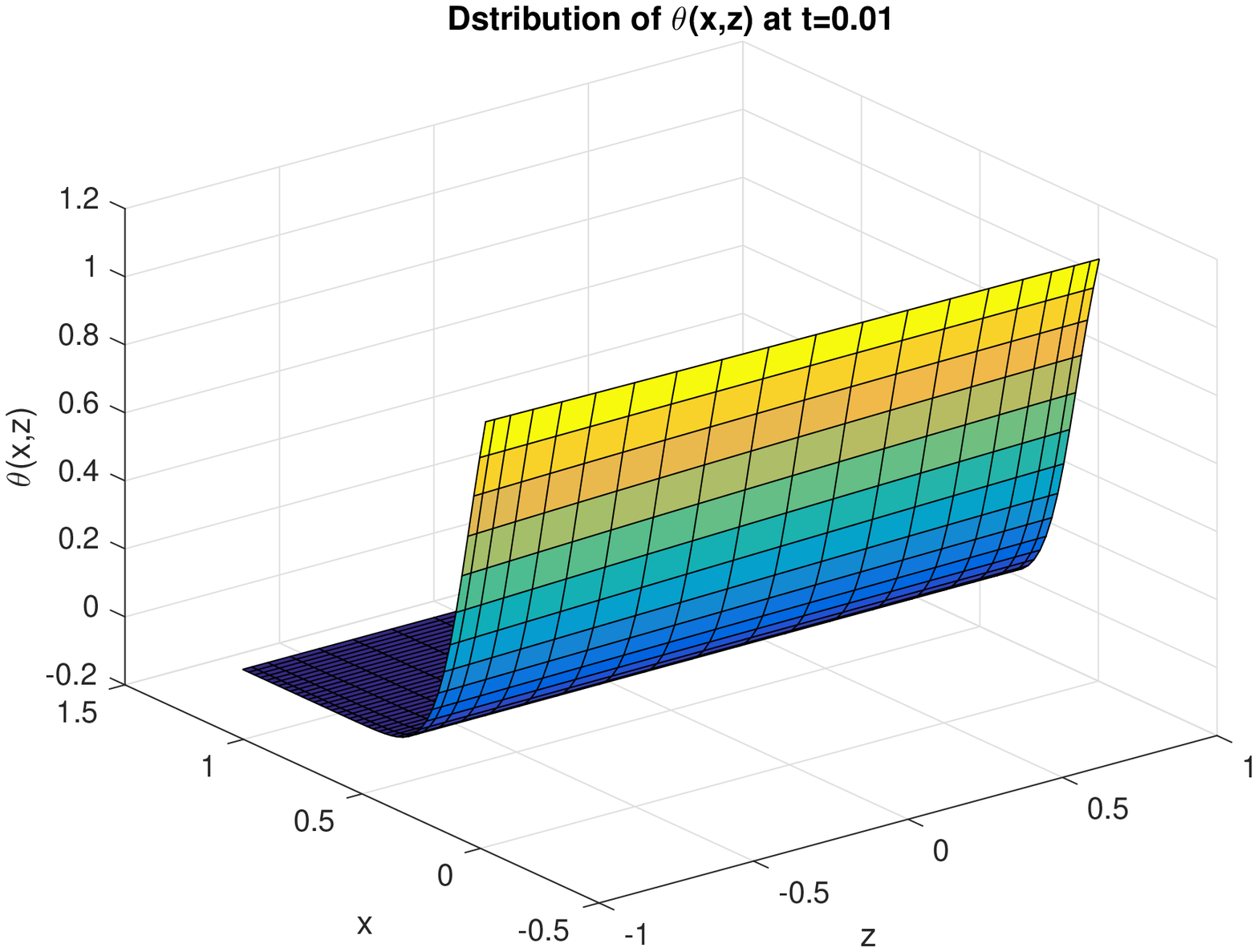}\includegraphics[scale=0.3]{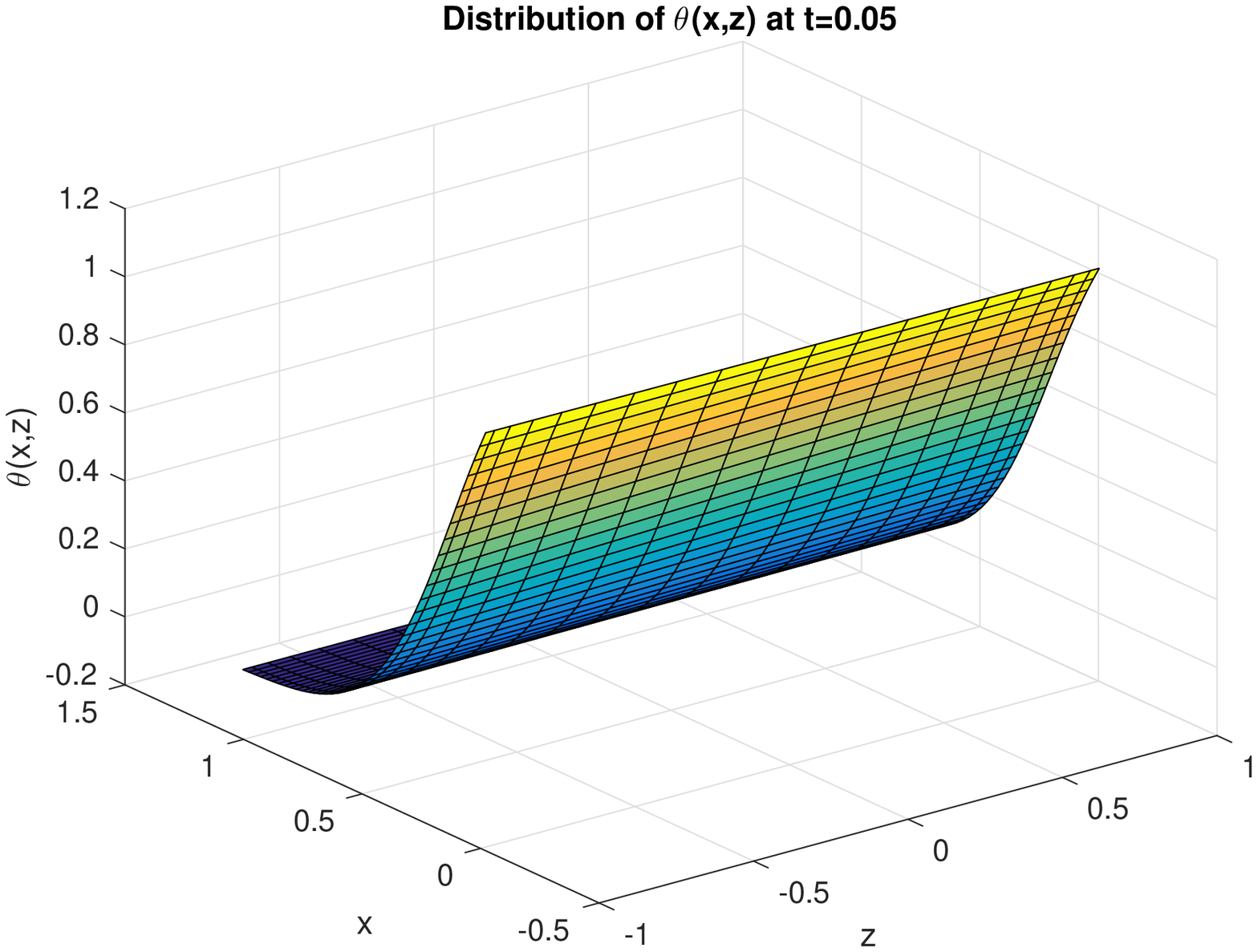}\includegraphics[scale=0.3]{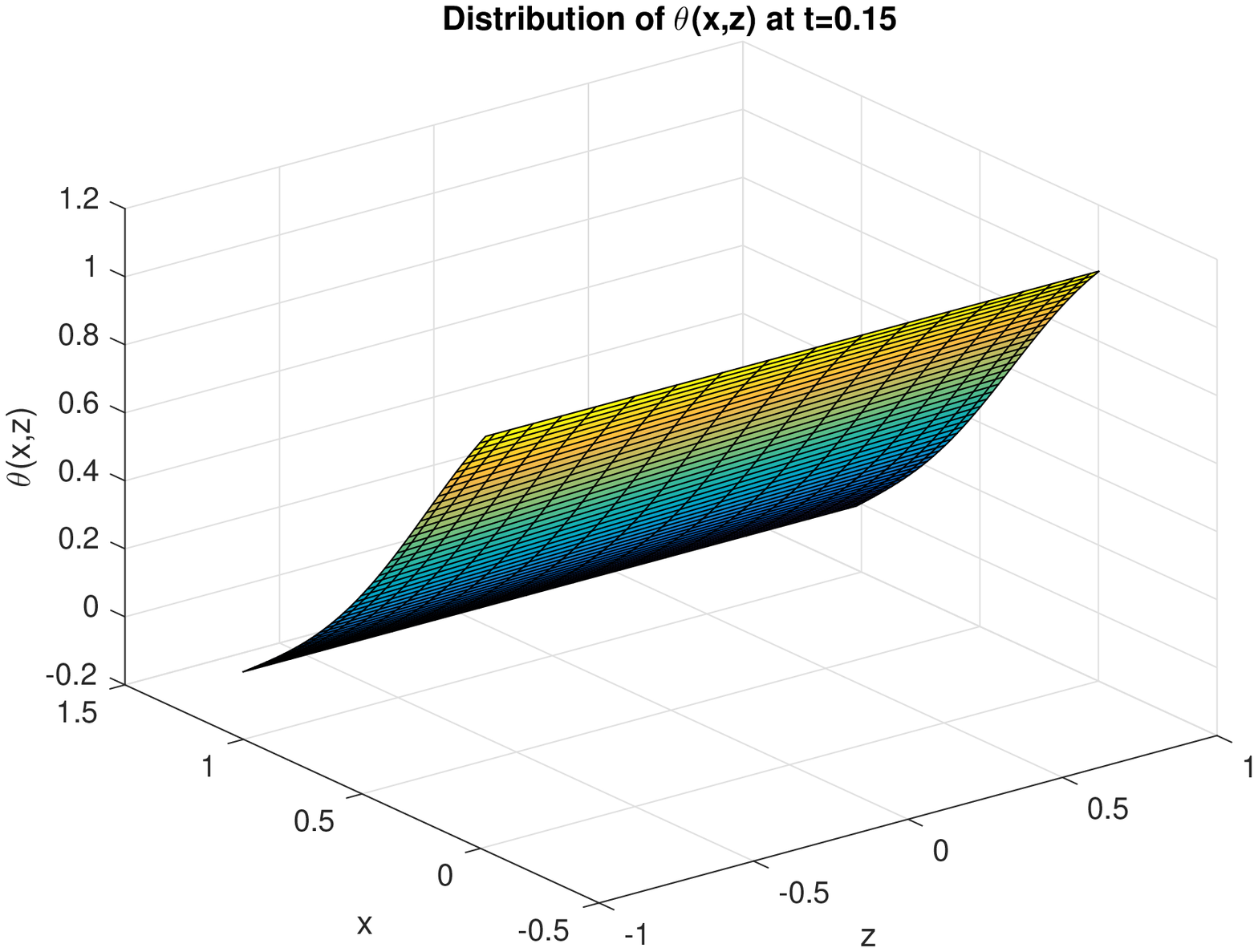}
\caption{Test 2. Distribution of the temperature $\theta(x,z)$ to the radiative heat transfer equation by gPC-SSP2 with $NT=11, 57, 171$.}
\end{figure}

\section{Conclusions}
In this article, for linear transport and radiative heat transfer equations with random inputs,  new generalized polynomial chaos based Asymptotic-Preserving stochastic Galerkin schemes are introduced. Compared with previous methods for these problems, our new method use the implicit-explicit (IMEX) time discretization to gain higher order accuracy, and by using a modified diffusion operator based penalty method, a more relaxed stability condition--a hyperbolic, rather than parabolic, CFL condition--is achieved when  mean free path is small, in the diffusive regime.  

These schemes allow efficient computation of random transport equations that contain both uncertainties and multiscales, allowing Knudsen number independent time step, mesh size and degree of polynomials, with a spectral accuracy in the random space.They can be used efficiently for all range of Knudsen numbers.

There remain many issues to be resolved. First, we only studied one dimensional problem in space and velocity.  The next step will be to extend the methods to high dimensions. In addition, multi-dimensional challenge is also presented in the random space.  Sparse grids could be used, as was done in \cite{ShuHuJin}, but it remains to be explored.

\bibliographystyle{plain}

\begin{thebibliography}{1}

\bibitem{BSS}
C. Bardos, R. Santos, and R. Sentis, Diffusion approximation and computation
of the critical size. Trans. Amer. Math. Soc., 284, 617--649,
1984.

\bibitem{boscarino2013implicit}
S. Boscarino, L. Pareschi, and G. Russo,
\newblock Implicit-explicit Runge--Kutta schemes for hyperbolic systems and
  kinetic equations in the diffusion limit.
\newblock  SIAM Journal on Scientific Computing, 35(1):A22--A51, 2013.

\bibitem{boscarino2016}
S. Boscarino, L. Pareschi, and G. Russo,
\newblock A unified IMEX Runge-Kutta approach for hyperbolic systems with multiscale relaxation.
\newblock  {\em Preprint}, 2017.

\bibitem{BGP}
F. Bouchut, F. Golse and M. Pulvirenti,
Kinetic Equations and Asymptotic Theory.
Seriers in Applied Math.,  Ed. B. Perthame and L. Desvillettes, Elsevier,
2000.

\bibitem{BC} C.~Buet, S.~Cordier, An asymptotic preserving scheme for hydrodynamics radiative transfer
models, Numerische Math. 108, 2, 199--221, 2007.

\bibitem{CGLV} {J.A.~Carrillo, T.~Goudon, P.~Lafitte, F.~Vecil}, Numerical schemes of diffusion
asymptotics and moment closures for kinetic equations. J. Sci.
Comput., 36, 113-149, 2008.

\bibitem{ChenLiuMu}
Z. Chen, L. Liu, and L. Mu, DG-IMEX Stochastic Galerkin schemes for Linear Transport Equation with Random Inputs and Diffusive Scalings. {\em Preprint}, 2016.

\bibitem{CZ} K.M.~Case, P.F.~Zweifel, {Linear Transport Theory.} Addison-Wesley, Reading, MA, 1997.

\bibitem{PN} F.~Cavalli, G.~Naldi, G.~Puppo, M.~Semplice, High order relaxation schemes for non linear diffusion problems.
SIAM Journal on Numerical Analysis, 45  2098-2119, 2007.

\bibitem{Cha} S.~Chandrasekhar, {Radiative Transport.} Dover, New York, 1960.

\bibitem{Deg1}
P. Degond, Asymptotic-preserving schemes for fluid models of plasmas. In Numerical models
for fusion, volume 39/40 of Panor. Syntheses, pages 1--€"90. Soc. Math. France, Paris, 2013.

\bibitem{Deg2}
P. Degond and F. Deluzet, Asymptotic-preserving methods and multiscale models for
plasma physics. J. Comp. Phys., 2017 (to appear).

\bibitem{Acta}
{G.~Dimarco, L.~Pareschi}, {Numerical methods for kinetic equations}, Acta Numerica 23, 369--520, 2014.

\bibitem{GT} {L.~Gosse, G.~Toscani}, {Asymptotic-preserving and well-balanced
schemes for radiative transfer and the Rosseland approximation}.
Numer. Math., 98, 2, 223--250, 2004.

\bibitem{HJ-Boltzmann}
J. Hu and S. Jin, A stochastic Galerkin method for the Boltzmann equation with uncertainty , J. Comp. Phys. 315, 150--168, 2016.

\bibitem{Jin_rev}
S. Jin, Asymptotic preserving (AP) schemes for multiscale kinetic and
hyperbolic equations: a review. Lecture Notes for Summer School on
 Methods and Models of Kinetic Theory €(M\&MKT), Porto Ercole (Grosseto,
Italy), 177--216, 2010.

\bibitem{JLM}
S. Jin, J.-G. Liu and Z. Ma, Uniform spectral convergence of the stochastic Galerkin method for the linear transport equations with random inputs in diffusive regime and a micro-macro decomposition based asymptotic preserving method. 
 Research in Math. Sci., to appear. 

\bibitem{JinLiu}
S. Jin and L. Liu,
 An Asymptotic-Preserving Stochastic Galerkin Method for the semiconductor Boltzmann equation with random inputs and diffusive scalings,
 Multiscale Model. Simul. 15, 157-183, 2017. 

\bibitem{JinLu}
S. Jin and H. Lu,
\newblock An asymptotic-preserving stochastic Galerkin method for the radiative
  heat transfer equations with random inputs and diffusive scalings.
J. Comp. Phys. 334, 182--206, 2017.


\bibitem{jin2000uniformly}
S. Jin, L. Pareschi, and G. Toscani,
\newblock Uniformly accurate diffusive relaxation schemes for multiscale
  transport equations.
\newblock  SIAM J. Num. Anal., 38(3),913--936, 2000.

\bibitem{jin2015asymptotic}
S. Jin, D. Xiu, and X. Zhu,
\newblock Asymptotic-preserving methods for hyperbolic and transport equations
  with random inputs and diffusive scalings.
\newblock J. Comp. Phys., 289, 35--52, 2015.

\bibitem{Klar} A. Klar,
An asymptotic-induced scheme for nonstationary transport
equations in the diffusive limit. SIAM J. Numer. Anal., 35, 1073--1094,
1998.

\bibitem{klar2001numerical}
A. Klar and C. Schmeiser,
\newblock Numerical passage from radiative heat transfer to nonlinear diffusion
  models.
\newblock  Math. Models Methods Appl. Sci.,
  11, 749--767, 2001.
  
\bibitem{LS} P.~Lafitte, G.~Samaey, Asymptotic-preserving projective integration schemes for kinetic equations in the diffusion limit, SIAM J. Sci. Comput., 34(2), A579--A602. 2012.

\bibitem{LK}
E.W. Larsen and J. B. Keller, Asymptotic solution of neutron
transport problems for small mean free paths. J. Math. Phys.,
15, 75€--81, 1974.

\bibitem{LMM1} 
E.W. Larsen and J. E. Morel, Asymptotic solutions of numerical
transport problems in optically thick, diffusive regimes. II. J. Comput.
Phys., 83, 212--236, 1989.

\bibitem{LMM2}
E.W. Larsen, J. E. Morel, and W.F. Miller, Jr., Asymptotic
solutions of numerical transport problems in optically thick, diffusive
regimes. J. Comput. Phys., 69, 283--324, 1987.


\bibitem{LM} M.~Lemou, L.~Mieussens,  A new asymptotic preserving scheme based on micro-macro formulation for linear kinetic equations in the diffusion limit. SIAM J. Sci. Comp. 31,  334--368, 2010.

\bibitem{LiWang}
Q. Li and L. Wang, Uniform regularity for linear kinetic equations with random input based on hypocoercivity, preprint. ArXiv: 1612.01219.

\bibitem{NP} G.~Naldi, L.~Pareschi, Numerical schemes for hyperbolic systems of conservation laws with stiff diffusive relaxation,  SIAM J.  Num. Anal.,
37,  1246--1270, 2000.

\bibitem{A} L.~Pareschi, G.~Russo, Implicit-Explicit Runge-Kutta schemes and applications to hyperbolic systems with relaxations. J. Sci. Comp. 25, 129--155, 2005.


\bibitem{ShuHuJin}
R. Shu, J. Hu and S. Jin,
\newblock A stochastic Galerkin method for the Boltzmann equation with high
  dimensional random inputs using sparse grids.
\newblock {\em Preprint}, 2016.

\bibitem{Spohn}
H. Spohn,  Large Scale Dynamics of Interacting Particles, Springer-Verlag, 
1991.

\bibitem{JiangXu}
W. Sun, S. Jiang, and K. Xu, An asymptotic preserving unified
gas kinetic scheme for gray radiative transfer equations. J. Comput. Phys.,
285, 265--279, 2015.

\bibitem{ZhuJin}
Y. Zhu and S. Jin, The Vlasov-Poisson-Fokker-Planck system with uncertainty and a one-dimensional asymptotic-preserving method. {\em Preprint}, 2016. 

\end{thebibliography}

\end{document}